\documentclass{amsproc}

\usepackage[arrow, curve, frame, matrix, graph]{xy}
\usepackage{amssymb}

\newtheorem{theorem}{Theorem}[section]
\newtheorem{lemma}[theorem]{Lemma}
\newtheorem{proposition}[theorem]{Proposition}
\newtheorem{corollary}[theorem]{Corollary}

\theoremstyle{definition}
\newtheorem{definition}[theorem]{Definition}
\newtheorem{example}[theorem]{Example}
\newtheorem{examples}[theorem]{Examples}
\newtheorem{non-example}[theorem]{Non-Example}

\theoremstyle{remark}
\newtheorem{remark}[theorem]{Remark}

\newcommand{\ca}[1]{\mathcal{#1}}

\newcommand{\Op}[2]{{\textnormal{Op}}(#1,#2)}
\newcommand{\N}{\mathbb{N}}

\newcommand{\id}{\textnormal{id}}

\newcommand{\tensor}{\otimes}
\newcommand{\iso}{\cong}
\newcommand{\catequiv}{\simeq}
\newcommand{\op}[1]{{#1}^{\textnormal{op}}}

\newcommand{\Shv}[1]{{\textnormal{Shv}}(#1)}
\newcommand{\Mod}[1]{{\textnormal{Mod}}(#1)}
\newcommand{\ArbMod}[2]{{\textnormal{Mod}}(#1,#2)}
\newcommand{\el}{\textnormal{el}}
\newcommand{\EL}{\textnormal{EL}}
\newcommand{\ladj}{\dashv}

\newcommand{\Cat}{\mathbf{Cat}}
\newcommand{\Set}{\mathbf{Set}}
\newcommand{\CAT}{\mathbf{CAT}}
\newcommand{\SET}{\mathbf{SET}}
\newcommand{\PSh}[1]{[{\op {#1}},{\Set}]}
\newcommand{\PSH}[1]{[{\op {#1}},{\SET}]}

\newcommand{\G}{\mathbb{G}}
\newcommand{\Glob}{\PSh \G}
\newcommand{\GLOBSET}{\PSH \G}
\newcommand{\GLOBCAT}{[{\op \G},{\CAT}]}
\newcommand{\CART}{\mathbf{PB}}
\newcommand{\TWOCAT}{\mathbf{2CAT}}
\newcommand{\Br}{\mathbf{Br}}
\newcommand{\Sym}{\mathbf{Sym}}
\newcommand{\TREE}{\ca T}
\newcommand{\Sol}{\textnormal{Sol}}
\newcommand{\PreOrd}{\mathbf{PreOrd}}
\newcommand{\Tr}{\mathbf{Tr}}

\newcommand{\wCat}{{\omega}{\textnormal{-}}{\Cat}}

\newcommand{\Simp}{\mathbf{\Delta}}
\newcommand{\solleq}{\blacktriangleleft}

\newcommand{\Span}{\mathbf{Span}}
\newcommand{\Coll}[1]{#1{\textnormal{-Coll}}}

\newcommand{\Alg}[1]{#1{\textnormal{-Alg}}}
\newcommand{\sAlg}[1]{#1{\textnormal{-Alg}_s}}
\newcommand{\PsAlg}[1]{{\textnormal{Ps-}}#1{\textnormal{-Alg}}}
\newcommand{\NormPsAlg}[1]{{\textnormal{Ps}}_{0}{\textnormal{-}}
#1{\textnormal{-Alg}}}

\newcommand{\PsMon}[1]{{\textnormal{PsMon(}}#1{\textnormal{)}}}
\newcommand{\Mon}[1]{{\textnormal{Mon(}}#1{\textnormal{)}}}

\begin{document}

\title{Operads within monoidal pseudo algebras}

\author{Mark Weber}
\address{Department of Mathematics and Statistics,
University of Ottawa}
\email{Mark.Weber@science.uottawa.ca}
\thanks{}


%
\bibliographystyle{amsalpha}
\maketitle
\begin{abstract}
A general notion of operad is given,
which includes:
\begin{enumerate}
\item  the operads that arose in algebraic topology in the 1970's
to characterise loop spaces.
\item  the higher operads of Michael Batanin \cite{Bat98}.
\item  braided and symmetric analogues of Batanin's operads
which are likely to be important in the study of
weakly symmetric higher dimensional monoidal categories.
\end{enumerate}
The framework of this paper,
links together $2$-dimensional monad theory,
operads, and higher dimensional algebra, in a natural way.
\end{abstract}

\section{Introduction}

Operads arose first in the early 1970's in algebraic topology
\cite{BV73} \cite{May72}
to keep track of the combinatorial data
that characterises infinite loop spaces.
In the most basic situation,
one has a braided monoidal category $(\ca V,I,\tensor)$,
and defines an operad to be a sequence of
objects $(p_n : n \in \N)$ of $\ca V$,
together with maps
\[ I \rightarrow p_1 \]
\[ p_k \tensor (p_{n_1} \tensor ... \tensor p_{n_k})
\rightarrow p_n \]
where $n = \sum_i n_i$.
This data satisfies some axioms, that 
ensure that it is sensible to regard each object
$p_n$, as an object of $n$-ary operations,
and the maps as expressing the process
of substitution of operations.
The $p_n$ and the corresponding maps
are called a \emph{non-symmetric operad
within a braided monoidal category}.
For applications, $\ca V$ can be some category of spaces,
chain complexes, differential graded algebras, or simplicial sets.
Typically, $\ca V$ is actually a symmetric monoidal category,
one has symmetric group actions
on each $p_n$, and asks that these actions
be compatible with the substitution.
Such an operad is known as a symmetric operad
within a symmetric monoidal category.

Beginning with insights of Todd Trimble \cite{Trb},
and then in the work of Michael Batanin \cite{Bat98},
operads were shown to be fundamental
for the explicit combinatorial description
of higher dimensional categorical structures.
However, the operads used in higher dimensional algebra
are typically somewhat more intricate than those originally
conceived to characterise loop spaces, although
the basic idea of formalising some notion of substitution
remains the same.

In \cite{Bat98}
as part of an approach to defining weak $\omega$-categories,
Batanin conceived of a notion of \emph{higher operad} internal to
a structure he called an augmented monoidal globular category.
This new operad notion is more complicated for two reasons.
First, monoidal categories are replaced by the more complicated
augmented monoidal globular categories.
Second, natural numbers $n$ as the place holders of the objects $p_n$
of the sequence,
are replaced by trees.
Just as addition of natural numbers may be regarded
as a consequence of the notion of monoidal category,
in that $\N$ with its addition is the strict monoidal
category freely generated by one object,
trees and their arithmetic operations (pasting of trees),
are encapsulated by the notion of monoidal globular category.

The notion of operad defined in this paper
formalises this phenomenon in the following way.
One begins with a $2$-monad $T$ on a $2$-category $\ca K$
whose job is two-fold:
\begin{enumerate}
\item  To describe the external structure within which
the corresponding operads live.
For example, to define non-symmetric operads
one takes $T$ to be the $2$-monad $\ca M$ on $\CAT$
whose algebras are monoidal categories.
Non-symmetric operads live
inside braided monoidal categories,
which are expressed here as \emph{monoidal pseudo algebras}
for the $2$-monad $\ca M$.
\label{key-idea}
\item  To encapsulate the ``indexing type''.
For example a sequence $p$ of objects in $\ca V$
is nothing but a functor $p:{\ca M}(1){\rightarrow}{\ca V}$,
because $\ca M(1) = \N$.
The monad structure of $\ca M$ expresses the addition
of natural numbers, which is necessary for the definition of
non-symmetric operad.
\end{enumerate}
An operad is then defined relative to $T$.
In this way, a unified formalism for the operads
originally considered in algebraic topology,
and those of interest to higher dimensional algebra,
is achieved, with different operad notions obtained
by varying $T$.

The idea central to this definition, is to regard
the external structure as composed of two parts:
a pseudo algebra structure for the $2$-monad,
together with a compatible pseudo monoid structure.
Taken together one has the notion of monoidal pseudo algebra
described in this paper.
The origin of this idea is in the observation
that when one describes the substitution maps for
a non-symmetric operad
\[ p_k \tensor (p_{n_1} \tensor ... \tensor p_{n_k})
\rightarrow p_n, \]
there are really two different types of tensor product
at work. One has a binary tensor product as in $p_k \tensor (...)$,
and $k$-ary tensor products as in $p_{n_1} \tensor ... \tensor p_{n_k}$.
The binary tensor product is formalised as the pseudo monoid structure, and
the $k$-ary tensor products are formalised as the pseudo $\ca M$-algebra structure.
Their compatibility implies that they can be identified
(as is the usual custom),
and that the resulting monoidal structure is braided.
The braiding is necessary for the expression of one of the operad axioms
(associativity of substitution).

The study of higher-dimensional braids
and tangles, as well as the homotopy
groups of spheres,
motivates the consideration of the various
notions of monoidal $n$-category.
A $k$-tuply monoidal weak $n$-category
is a weak $(n+k)$-category with one cell
in each dimension less than $k$.
Such a structure is considered as being
$n$-dimensional by reindexing appropriately,
that is,
by regarding the $m$-cells (for $m \geq k$)
of the original weak $(n+k)$-category
as $(m-k)$-cells in this new structure.
For example,
a $2$-tuply monoidal tricategory
is a braided monoidal category,
and a braiding is a subtler notion
of symmetry for monoidal categories than the usual one.
Thus, one expects these $k$-tuply monoidal weak $n$-categories
in general, to be higher dimensional monoidal categories
which possess still more subtle symmetry.

Motivated by insights from homotopy theory,
\cite{BD98b} give three hypotheses
relating such structures to quantum topology.
The first is that
the $n$-cells of the free
$k$-tuply monoidal weak $n$-category
on one object correspond to
``$n$-braids in $(n+k)$ dimensions'',
which are certain $n$-dimensional surfaces
embedded in the $(n+k)$-dimensional cube.
The case $n=1$ and $k=2$ gives the usual
definition of braid, which corresponds
to the morphisms of the free braided monoidal category
on one object.
Second is the
corresponding hypothesis for weak $n$-groupoids,
which relates to the fundamental $n$-groupoid
of the $k$-fold loop space of the $k$-sphere.
Finally, it is predicted that
the $n$-cells of
the free $k$-tuply monoidal weak $n$-category
\emph{with duals}
correspond to ``framed $n$-tangles in $(n+k)$ dimensions'',
which again are certain $n$-dimensional surfaces
embedded in the $(n+k)$-dimensional cube.
This time the case $n=1$ and $k=2$
corresponds to the usual definition of tangle
which has been shown to correspond to the free
braided monoidal category with duals on one object.

An important motivation for this work
is to define braided and symmetric analogues
of higher operads,
to facilitate the study of these weakly symmetric
higher dimensional categories.
With the general operad definition at our disposal,
this problem is reduced to finding appropriate
$2$-monads,
which blend together the combinatorics of higher operads
with braids and symmetries in a natural way.
The $2$-monad that parametrises Batanin's higher operads
is denoted by $\ca T$ and acts on the $2$-category
$\GLOBCAT$ of globular categories.
Moreover there are $2$-monads $\ca B$ and $\ca S$
on $\CAT$
which parametrise braided and symmetric operads in the usual sense.

The appropriate $2$-monads alluded to above
are obtained by regarding $\ca B$ and $\ca S$
as $2$-monads on $\GLOBCAT$ in an obvious way,
and seeing that there are distributive
laws between between these $2$-monads and $\ca T$.
The existence of these distributive laws
is deduced from an alternative description
of the category $\omega$-$\Cat$,
of strict $\omega$-categories, due to Clemens Berger \cite{Ber00}.

As observed in \cite{Lei} and \cite{Str00},
the higher operads which are actually used
in \cite{Bat98} to define weak $\omega$-categories,
all live in a particular augmented monoidal globular category
called $\Span$, and admit a far simpler description.
One has the monad on $\ca T_0$ on $\Glob$ the category of globular sets,
whose algebras
are strict $\omega$-categories, and a higher operad in $\Span$
amounts to a cartesian monad morphism $\phi_0:R_0{\rightarrow}{\ca T_0}$.
That is, a monad $R_0$ on $\Glob$, and
a natural transformation $\phi_0$ which is compatible with the monad structures,
and whose naturality squares are pullbacks.
We call such higher operads \emph{basic higher operads}.
On the other hand, \cite{Bat02}
uses the full generality of the higher operad notion
for applications to loop spaces.
So, while basic higher operads suffice
for the definition of weak $\omega$-category presented in \cite{Bat98},
it seems that general higher operads are important for
applications.

In this paper we must speak of two set theoretic universes
$\ca U_1 \in \ca U_2$
and distinguish between
$\Set$, the category of $\ca U_1$ small sets and functions between them,
and $\SET$, the category of $\ca U_2$ small sets.
Similarly we distinguish between the corresponding $2$-categories
$\Cat$ and $\CAT$ of categories.
So $\Set$ and $\Cat$ may be regarded as objects of $\CAT$,
as may many of the other categories that one
encounters in applications of operads: categories of spaces,
chain complexes, differential graded algebras and simplicial sets.
The reason for this distinction is that the $2$-monads $T$
on $2$-categories $\ca K$ that parametrise operad notions,
apply at the $\ca U_2$ level. For example,
the monad $\ca M$ on $\CAT$ that parametrises non-symmetric operads
within braided monoidal categories.
The braided monoidal categories within which our operads
live are objects of $\CAT$.

Having made this distinction, it is worth noting on the other hand,
that most general categorical and combinatorial constructions do not
depend on such considerations, that is, they are ``universe-insensitive''.
This is part of the reason why such size issues
are often glossed over. However in \cite{Str00}, such distinctions
are shown to be pertinent to the organisation of the combinatorics
and category theory which underlies higher dimensional algebra.
As part of such distinctions,
we have used the notation $\phi_0:R_0{\rightarrow}{\ca T_0}$
to denote a basic higher operad,
which is a morphism of monads on $\Glob$.
Associated to $\phi_0$ is the morphism of $2$-monads
$\phi:R{\rightarrow}{\ca T}$ on $\GLOBCAT$,
obtained from $\phi_0$ by changing universes
and taking category objects.
This notation is convenient for us, because
the distributivity of braids and symmetries with $\ca T$
mentioned above,
is actually more general -- one can replace $\ca T$ with $R$
for any basic higher operad $\phi_0$.
For example, $R_0$ could be a monad on $\Glob$ whose algebras
are \emph{weak} $\omega$-categories in the sense of \cite{Bat98}.
In this way one also has weakened versions of
higher operads, as well as their braided and
symmetric analogues, captured by our general formalism.

This paper is organised as follows.
Sections (\ref{2-monad}) and (\ref{ps-monoid})
review $2$-monads and their algebras,
and pseudo monoids,
assuming familiarity with the usual categorical notions
of monad and monoid.
Monoidal pseudo algebras are introduced in section (\ref{mon-psalg}),
and in section (\ref{operads}),
operads and their algebras are defined in full generality.
The examples presented in sections (\ref{2-monad})--({\ref{operads}),
taken together, exhibit how the conventional operad notions
are captured by our formalism.
Then in section (\ref{sec:higher-op}),
after briefly recalling the relevant background on the
globular approach to higher category theory,
the higher operads of \cite{Bat98}
are described as instances of our general operads.
We begin section ({\ref{higher-sym-ops})
by recalling the characterisation
from \cite{Ber00}
of the category algebras of a basic higher operad.
This is then re-expressed in the language of sketches,
which then allows the easy explanation
of the formal distributivity of symmetries and braids with
basic higher operads.

\section{$2$-monads and pseudo algebras}\label{2-monad}

Recall that a $2$-monad $(T,\eta,\mu)$ on a $2$-category $\ca K$
consists of an endo-$2$-functor $T$ of $\ca K$,
together with $2$-natural transformations $\eta:1{\implies}T$
and ${\mu}:T^2{\implies}T$, called the unit and multiplication, so that
\[ \begin{xy}
(0,0);<2em,0em>:<0em,2em>::
(-1,0)*+!R{\xybox{
\xymatrix{T \ar[r]^-{\eta{T}} \ar[dr]_{1_T} & {T^2} \ar[d]|{\mu}
& T \ar[l]_-{T\eta} \ar[dl]^{1_T} \\ & T}}},
(1,0)*+!L{\xybox{
\xymatrix{{T^3} \ar[r]^-{\mu{T}} \ar[d]_{T\mu} & {T^2} \ar[d]^{\mu} \\
{T^2} \ar[r]_-{\mu} & T}}}
\end{xy} \]
commute.
We shall allow the usual abuse of referring to the $2$-monad $T$,
omitting reference to the unit and multiplication.
Most of the examples of $2$-monads of interest to us shall now be described,
and for many more examples, the reader may consult \cite{BKP}.
\begin{examples}\label{Monad-ex}
\begin{enumerate}
\item  Every monad $(T,\eta, \mu)$ on a category $\ca E$
can be regarded as a $2$-monad, by regarding $\ca E$
as a locally discrete $2$-category (that is, one with only identity $2$-cells).
\label{loc-discrete}
\item  For each $\ca K$ one obtains the identity monad $1_{\ca K}$ on $\ca K$,
by taking $T$, $\eta$, and $\mu$ to be identities.
\label{id}
\item  Let $\ca E$ be a category with pullbacks,
and $(T,\eta,\mu)$ be a monad on $\ca E$ such that $T$ preserves pullbacks.
One can then take $\ca K$ to be the $2$-category
${\Cat}(\ca E)$ of categories internal to $\ca E$.
This process of taking the $2$-category of internal categories
is the object map of a $2$-functor
\[ \xymatrix{{\CART} \ar[rr]^-{\Cat} && {\TWOCAT}} \]
from the $2$-category of categories with pullbacks,
pullback preserving functors and natural transformations between them,
to the $2$-category of $2$-categories, $2$-functors and $2$-natural
transformations.
Applying this $2$-functor, one obtains a $2$-monad $\Cat(T)$
on $\Cat(\ca E)$.
\label{Cat-of-monad}
\item  As an instance of (\ref{Cat-of-monad}),
take $\ca E$ to be $\SET$ and $T$ the monoid monad on $\SET$.
We denote by $\ca M$ the $2$-monad $\Cat(T)$ on $\CAT$.
An object of ${\ca M}(X)$ is a sequence of objects from $X$,
that is, a functor $x:n{\rightarrow}X$ where $n \in \N$ is being regarded
as the discrete category whose object set is $n = \{0,...,n-1\}$.
A morphism $f:x{\rightarrow}y$ in ${\ca M}(X)$ is a $2$-cell
\[ \xymatrix{n \ar@/^1pc/[r]^-{x} \save \POS?="dom" \restore
\ar@/_1pc/[r]_-{y} \save \POS?="cod" \restore & X
\POS "dom"; "cod" **@{};
? /-.5em/ \ar@{=>}^{f} ? /.5em/ }, \]
and so is just a sequence of maps in $X$.
The $1$ and $2$-cell mappings for $\ca M$
are obtained by composition in the evident fashion.
The unit for the monad picks out the sequences of length one,
and the multiplication is given by concatenation of sequences.
\label{Mon-on-Cat}
\item  Denote by ${\Br}_n$ the $n$-th braid group.
We shall denote by $\ca B$ the following $2$-monad
on $\CAT$.
An object of ${\ca B}(X)$ is again a sequence of objects of $X$.
A morphism between two sequences $x$ and $y$ of the same length $n$
consists of a braid on $n$-strings, whose strings are labelled
by arrows in $X$. More precisely, such a morphism consists of
$\beta \in {\Br}_n$, together with a $2$-cell
\[ \xymatrix{n \ar[rr]^-{\overline{\beta}} \ar[dr]_{x}
\save \POS?="dom" \restore && n \ar[dl]^{y}
\save \POS?="cod" \restore \\ & X
\POS "dom"; "cod" **@{};
? /-.5em/ \ar@{=>}^{f} ? /.5em/ }, \]
where $\overline{\beta}$ is the underlying permutation of $\beta$
regarded as a functor between discrete categories.
The $2$-functoriality of $\ca B$ and the unit work as with $\ca M$.
The multiplication is described by concatenation of sequences,
and substitution of braids into braids in the evident way.
\label{BrMon-on-Cat}
\item  Denote by ${\Sym}_n$ the $n$-th symmetric group.
We shall denote by $\ca S$ the following $2$-monad
on $\CAT$.
An object of ${\ca S}(X)$ is again a sequence of objects of $X$.
A morphism between two sequences $x$ and $y$ of the same length $n$
consists of a permutation on $n$-strings whose strings are labelled
by arrows in $X$. More precisely, such a morphism consists of
$\beta \in {\Sym}_n$, together with a $2$-cell
\[ \xymatrix{n \ar[rr]^-{\beta} \ar[dr]_{x}
\save \POS?="dom" \restore && n \ar[dl]^{y}
\save \POS?="cod" \restore \\ & X
\POS "dom"; "cod" **@{};
? /-.5em/ \ar@{=>}^{f} ? /.5em/ }, \]
where $\beta$ is being
regarded as a functor between discrete categories.
The $2$-functoriality of $\ca S$ and the unit work as with $\ca M$.
The multiplication is described by concatenation of sequences,
and substitution of permutations into permutations in the evident way.
\label{SymMon-on-Cat}
\end{enumerate}
\end{examples}
The important difference between
$2$-monads and ordinary monads on categories,
is that there are various weaker notions of algebra
in addition to the usual (Eilenberg-Moore) algebras
for a monad. This makes $2$-monad theory
a natural choice of formalism when one wishes to consider
coherently defined categorical structures.
In this work we shall consider pseudo algebras
and pseudo morphisms -- where one replaces equality
between composite arrows
in the axiomatic definition of the objects and arrows of
$\Alg T$, the category of Eilenberg-Moore algebras for $T$, by isomorphisms.
\begin{definition}
Let $(T,\eta,\mu)$ be a $2$-monad on a $2$-category $\ca K$.
A \emph{pseudo $T$-algebra structure}
$(a,\alpha_0,\alpha)$
on an object $A \in \ca K$ consists of
a $1$-cell
$a:{T}A{\rightarrow}A$
and invertible $2$-cells
\[ \begin{xy}
(0,0);<2em,0em>:<0em,2em>::
(-1,0)*+!R{\xybox{
\xymatrix{{T^2A}
\ar[r]^-{{\mu}_{A}} \ar[d]_{{T}a}
\save \POS?="domass" \restore
& {{T}A} \ar[d]^{a}
\save \POS?="codass" \restore \\
{{T}A} \ar[r]_-{a}
& A
\save \POS"domass";"codass"**@{}
?  /-.5em/ \ar@{=>}^-{\alpha} ?  /.5em/ \restore
}}},
(1,0)*+!L{\xybox{
\xymatrix{A \ar[rr]^-{{\eta}_{A}} \ar[dr]_{1_A}
\save \POS?="codlam" \restore
&& {{T}A} \ar[dl]^{a}
\save \POS?="domboth" \restore \\
& A
\POS "codlam";"domboth"**@{}
?  /-.5em/ \ar@{=>}^-{\alpha_0} ?  /.5em/ \restore
}}},
\end{xy} \]
in $\ca K$ satisfying
\[ \begin{xy}
(0,0);<2em,0em>:<0em,2em>::
(-1,0)*-!R{\xybox{
\xymatrix @C=1em {
& {{T}^{2}A} \ar[rr]^-{{\mu}_{A}}
&& {{T}A} \ar[dr]^-{a}
\save \POS?="m" \restore
\\
{{T}^{3}A} \ar[rr]
\ar[ur]^-{{\mu}_{{T}A}}
\save \POS?="11m" \restore
\ar[dr]_-{{T}^{2}a}
\save \POS?="m11" \restore
&& {{T}^{2}A}
\ar[ur] \save \POS?="1m" \restore
\ar[dr]|{{T}a} \save \POS?="m1" \restore
&& {A} \\
& {{T}^{2}A} \ar[rr]_-{{T}a}
&& {{T}A} \ar[ur]_-{a}
\POS "1m"; "11m" **@{}
? /-.35em/ \ar@{=} ? /.35em/
\POS "m11"; "m1" **@{}
? /-.5em/ \ar@{=>}^{{T}{\alpha}} ? /.5em/
\POS "m1"; "m" **@{}
? /-.5em/ \ar@{=>}^{\alpha} ? /.5em/
}}}="lhs";
(1,0)*-!L{\xybox{
\xymatrix @C=1em {
& {{T}^{2}A} \ar[rr]^-{{\mu}_{A}}
\ar[dr] \save \POS?="m1" \restore
&& {{T}A} \ar[dr]^-{a}
\save \POS?="mtop" \restore \\
{{T}^{3}A}
\ar[ur]^-{{\mu}_{{T}A}}
\ar[dr]_-{{T}^{2}a}
\save \POS?="m11" \restore
&& {{T}A} \ar[rr]^-{a}
&& {A} \\
& {{T}^{2}A} \ar[rr]_-{{T}a}
\ar[ur] \save \POS?="1m" \restore
&& {{T}A} \ar[ur]_-{a}
\save \POS?="mbot" \restore
\POS "m11"; "m1" **@{}
? /-.35em/ \ar@{=} ? /.35em/
\POS "m1"; "mtop" **@{}
? /-.5em/ \ar@{=>}^{\alpha} ? /.5em/
\POS "mbot"; "1m" **@{}
? /-.5em/ \ar@{=>}_{\alpha} ? /.5em/
}}}="rhs"
**@{} ? /-.5em/ \ar@{=} ? /.5em/
\end{xy} \]
and
\[ \begin{xy}
(0,0);<2em,0em>:<0em,2em>::
(-1,0)*+!R{\xybox{
\xymatrix{{TA} \ar[dd]_{1_{TA}} \save \POS?="domTalp0" \restore
\ar[dr] \save \POS?="domeq" \restore \ar[rr]^-{1_{TA}}
&& {TA} \ar[dd]^{a} \save \POS?(.75)="codalp" \restore \\
& {T^2A} \save \POS="codTalp0" \restore
\ar[ur]|{\mu_A} \save \POS?="codeq" \restore
\ar[dl]|{Ta} \save \POS?="domalp" \restore \\
TA \ar[rr]_-{a} && A
\POS "domeq";"codeq" **@{} ? *{=}
\POS "domTalp0";"codTalp0" **@{}
? /-.75em/ \ar@{=>}^{T\alpha_0} ? /.25em/
\POS "domalp";"codalp" **@{}
? /-.5em/ \ar@{=>}^{\alpha} ? /.5em/}}}="lhs";
(1,0)*+!L{\xybox{*{1_a.}}}="rhs"
**@{} ? /-.5em/ \ar@{=} ? /.5em/
\end{xy} \]
The triple $(A,\alpha_0,\alpha)$
is referred to as a \emph{pseudo $T$-algebra.}
When $\alpha_0$ is an identity
the pseudo algebra is said to be \emph{normal}.
When in addition $\alpha$ is an identity,
we refind the usual notion of $T$-algebra, and
the algebra is said to be \emph{strict}.
\end{definition}
\begin{definition}
Let $(A,\alpha_0,\alpha)$ and $(A',{\alpha}'_0,{\alpha}')$
be pseudo $T$-algebras.
A \emph{strong $T$-morphism structure} for a $1$-cell
$f:A{\rightarrow}A'$ is an invertible $2$-cell
\[ \xymatrix{{{T}A} \ar[r]^-{a} \ar[d]_{{T}f}
\save \POS?="domfbar" \restore
& {A} \ar[d]^{f}
\save \POS?="codfbar" \restore \\
{{T}A'} \ar[r]_-{a'} & {A'}
\POS"domfbar"; "codfbar" **@{}
? /-.5em/ \ar@{=>}^{\overline{f}} ? /.5em/
} \]
satisfying
\[ \begin{xy}
(0,0);<2em,0em>:<0em,2em>::
(-1,0)*-!R{\xybox{
\xymatrix @C=1em {
& {{T}A} \ar[rr]^-{a}
&& {A} \ar[dr]^-{f}
\save \POS?="m" \restore \\
{{T}^{2}A} \ar[rr]^-{{T}a}
\ar[ur]^-{{\mu}_{A}}
\save \POS?="11m" \restore
\ar[dr]_-{{T}^{2}f}
\save \POS?="m11" \restore
&& {{T}A}
\ar[ur]|{a} \save \POS?="1m" \restore
\ar[dr]|{{T}f} \save \POS?="m1" \restore
&& {A'} \\
& {{T}^{2}A'} \ar[rr]_-{{T}a'}
&& {{T}A'} \ar[ur]_-{a'}
\POS "1m"; "11m" **@{}
? /-.5em/ \ar@{=>}_{\alpha} ? /.5em/
\POS "m11"; "m1" **@{}
? /-.5em/ \ar@{=>}^{{T}\overline{f}} ? /.5em/
\POS "m1"; "m" **@{}
? /-.5em/ \ar@{=>}^{\overline{f}} ? /.5em/
}}}="lhs";
(1,0)*-!L{\xybox{
\xymatrix @C=1em {
& {{T}A} \ar[rr]^-{a}
\ar[dr]|{{T}f} \save \POS?="m1" \restore
&& {A} \ar[dr]^-{f}
\save \POS?="mtop" \restore \\
{T^2A}
\ar[ur]^-{{\mu}_{A}}
\ar[dr]_-{{T}^{2}f}
\save \POS?="m11" \restore
&& {{T}A'} \ar[rr]^-{a'}
&& {A'} \\
& {{T}^{2}A'} \ar[rr]_-{{T}a'}
\ar[ur]|{{\mu}_{A'}} \save \POS?="1m" \restore
&& {{T}A'} \ar[ur]_-{a'}
\save \POS?="mbot" \restore
\POS "m11"; "m1" **@{}
? /-.35em/ \ar@{=} ? /.35em/
\POS "m1"; "mtop" **@{}
? /-.5em/ \ar@{=>}^{\overline{f}} ? /.5em/
\POS "mbot"; "1m" **@{}
? /-.5em/ \ar@{=>}_{{\alpha}'} ? /.5em/
}}}="rhs"
**@{} ? /-.5em/ \ar@{=} ? /.5em/
\end{xy} \]
and
\[ \begin{xy}
(0,0);<2em,0em>:<0em,2em>::
(-1,0)*-!R{\xybox{
\xymatrix @R=1.2em {& {TA} \ar[dr]^{a} \ar[dd]
\save \POS?="domfbar" \restore \save \POS?="codeq" \restore \\
A \ar[ur]^{\eta_A} \ar[dd]_{f} \save \POS?="domeq" \restore
&& A \ar[dd]^{f} \save \POS?="codfbar" \restore \\
& {TA'} \save \POS="codalpp0" \restore \ar[dr] \\
{A'} \ar[ur] \ar[rr]_-{1_{A'}} \save \POS?="domalpp0" \restore && {A'}
\POS "domfbar"; "codfbar" **@{}
? /-.5em/ \ar@{=>}^{\overline{f}} ? /.5em/
\POS "domalpp0"; "codalpp0" **@{}
? /-.6em/ \ar@{=>}^{{\alpha}'_0} ? /.2em/
\POS "domeq"; "codeq" **@{}
? /-.35em/ \ar@{=} ? /.35em/ }}}="lhs";
(1,0)*-!L{\xybox{
\xymatrix{A \ar[rr]^-{\eta_A} \ar[dr]_{1_A} \save \POS?="domalp0" \restore
&& {TA} \ar[dl]^{a} \save \POS?="codalp0" \restore \\
& A \ar[d]^{f} \\ & {A'}
\POS "domalp0"; "codalp0" **@{}
? /-.5em/ \ar@{=>}^{\alpha_0} ? /.5em/ }}}="rhs"
**@{} ? /-.5em/ \ar@{=} ? /.5em/
\end{xy} \]
The pair $(f,\overline{f})$ is called a \emph{strong $T$-morphism.}
We shall allow the notational abuse of referring to the ``strong $T$-morphism $f$'',
omitting any reference to $\overline{f}$.
When $\overline{f}$ is an identity, we refind the usual notion
of $T$-algebra morphism, and the $T$-morphism is said to be \emph{strict}
in this case.
\end{definition}
\begin{definition}
Let $f$ and $f'$ be strong $T$-morphisms
$(a,\alpha_0,\alpha)\rightarrow(a',{\alpha}'_0,{\alpha}')$.
A $2$-cell $\psi:f{\implies}f'$ is an \emph{algebra $2$-cell} when
\[ \begin{xy}
(0,0);<2em,0em>:<0em,2em>::
(-1,0)*-!R{\xybox{
\xymatrix{{{T}A} \ar@/_1pc/[d]_{{T}f}
\save \POS?="dompsi" \restore
\ar@/^1pc/[d]^{{T}f'}
\save \POS?="codpsi" \restore
\save \POS?="domphip2" \restore
\ar[r]^-{a}
& {A}
\ar@/^1pc/[d]^{f'}
\save \POS?="codphip2" \restore \\
{{T}A'} \ar[r]_-{a'} & {A'}
\POS "dompsi"; "codpsi" **@{}
? /-.5em/ \ar@{=>}^{T{\psi}} ? /.5em/
\POS "domphip2"; "codphip2" **@{}
? /-.5em/ \ar@{=>}^{\overline{f'}} ? /.5em/
}
}}="lhs";
(1,0)*-!L{\xybox{
\xymatrix{{{T}A} \ar@/_1pc/[d]_{{T}f}
\save \POS?="domphi2" \restore
\ar[r]^-{a}
& {A}
\ar@/_1pc/[d]_{f}
\save \POS?="dompsi2" \restore
\save \POS?="codphi2" \restore
\ar@/^1pc/[d]^{f'}
\save \POS?="codpsi2" \restore \\
{{T}A'} \ar[r]_-{a'} & {A'}
\POS "domphi2"; "codphi2" **@{}
? /-.5em/ \ar@{=>}^{\overline{f}} ? /.5em/
\POS "dompsi2"; "codpsi2" **@{}
? /-.5em/ \ar@{=>}^{\psi} ? /.5em/
}
}}="rhs"
**@{} ? /-.5em/ \ar@{=} ? /.5em/
\end{xy} \]
\end{definition}
With the evident compositions,
one defines the $2$-category $\PsAlg T$
to consist of pseudo $T$-algebras,
strong $T$-morphisms and algebra $2$-cells.
The full sub-$2$-category of $\PsAlg T$ consisting of
the normal pseudo algebras is denoted $\NormPsAlg T$.
The locally full sub-$2$-category of $\PsAlg T$ consisting of
the strict algebras and strict morphisms is denoted $\sAlg T$.
\begin{examples}
\begin{enumerate}
\item  The $2$-categories of strict and pseudo algebras coincide
for (\ref{Monad-ex})({\ref{loc-discrete}}),
being just the usual category of algebras for $T$
regarded as a locally discrete $2$-category.
\item  For any $\ca K$,
a strict algebra structure for $1_{\ca K}$ is vacuous.
A normal pseudo algebra structure is also vacuous.
A pseudo algebra structure on $X \in \ca K$ amounts
to $t:X{\rightarrow}X$
together with an isomorphism $t {\iso} 1_X$.
\item  The $2$-category of strict algebras for (\ref{Monad-ex})(\ref{Cat-of-monad})
is just $\Cat(\Alg T)$.
\item  A strict $\ca M$-algebra structure on a category $X$
is a strict monoidal structure.
A pseudo $\ca M$-algebra structure on a category $X$
is a monoidal structure, described in an \emph{unbiased} fashion.
That is, one supplies an $n$-ary tensor product for $n \in \N$,
and associated coherence isomorphisms.
For a normal pseudo $\ca M$-algebra structure,
the $1$-ary tensor product of $x \in X$ is $x$, rather than just
isomorphic to $x$.
There are various monoidal coherence results in the literature,
for example in \cite{Power}, \cite{LkCodesc}, \cite{Her00}
and \cite{HerCoh},
which are expressed in the language of pseudo algebras,
and so apply to many other situations.
In all these results, the inclusion $2$-functor
$\sAlg {\ca M}{\rightarrow}{\PsAlg {\ca M}}$
is seen to be a biequivalence.
In addition to these results,
one can also exhibit directly, a $2$-equivalence between
$\NormPsAlg {\ca M}$ and the $2$-category $\PsMon{\CAT}$
consisting of monoidal categories defined
in the usual (biased) way, by giving a binary tensor product
and a unit object.
Under this $2$-equivalence,
strong $\ca M$-morphisms coincide with the tensor functors of \cite{JS93},
and have been called strong monoidal functors elsewhere.
\item  A strict $\ca B$-algebra structure on a category $X$
is a braided strict monoidal structure, that is, a braided tensor category
in the sense of \cite{JS93} whose underlying monoidal category is strict.
A pseudo $\ca B$-algebra structure on a category $X$
amounts to a braided monoidal structure on $X$.
\item  Similarly, strict and pseudo algebras for $\ca S$
are symmetric strict monoidal categories and
symmetric monoidal categories respectively.
\end{enumerate}
\end{examples}
%

\section{Pseudo monoids}\label{ps-monoid}

Having described a well-known ``categorification'' of monad algebra,
we shall now consider one for the notion of monoid
in a monoidal category.
For us, it suffices to consider pseudo monoids within
$2$-categories with cartesian products in the $\CAT$-enriched sense,
rather than internal to a more general monoidal $2$-category.
Later, when we describe monoidal pseudo algebras
and the operads they contain,
this specialisation to $2$-categories with \emph{cartesian} products
becomes crucial.
The reason for this, as we shall see,
is that pseudo monoids within such $2$-categories
can be described \emph{representably}.
For the remainder of this section,
$\ca K$ is a $2$-category with finite products.
\begin{definition}
A \emph{pseudo monoid structure}
$(i,m,\alpha,\lambda,\rho)$
on $A \in \ca K$
consists of $1$-cells
\[ \xymatrix{{1} \ar[r]^-{i}
& A & {A{\times}A} \ar[l]_-{m}} \]
and invertible $2$-cells
\[ \begin{xy}
(0,0);<2em,0em>:<0em,2em>::
(-1,0)*+!R{\xybox{
\xymatrix{{A{\times}A{\times}A}
\ar[r]^-{1{\times}m} \ar[d]_{m{\times}1}
\save \POS?="domass" \restore
& {A{\times}A} \ar[d]^{m}
\save \POS?="codass" \restore \\
{A{\times}A} \ar[r]_-{m}
& A
\save \POS"domass";"codass"**@{}
?  /-.5em/ \ar@{=>}^-{\alpha} ?  /.5em/ \restore
}}},
(1,0)*+!L{\xybox{
\xymatrix{A \ar[r]^-{i{\times}1} \ar[dr]_{1}
\save \POS?="codlam" \restore
& {A{\times}A} \ar[d]|{m}
\save \POS?="domboth" \restore
& A \ar[l]_-{1{\times}i} \ar[dl]^{1}
\save \POS?="codrho" \restore \\
& A
\save \POS"domboth";"codlam"**@{}
?  /-.5em/ \ar@{=>}_-{\lambda} ?  /.5em/ \restore
\save \POS"domboth";"codrho"**@{}
?  /-.5em/ \ar@{=>}^-{\rho} ?  /.5em/ \restore
}}},
\end{xy} \]
in $\ca K$
satisfying the following two axioms:
\[ \begin{xy}
(0,0);<2em,0em>:<0em,2em>::
(-1,0)*-!R{\xybox{
\xymatrix @C=1em {
& {A^{3}} \ar[rr]^-{1{\times}m}
&& {A^{2}} \ar[dr]^-{m}
\save \POS?="m" \restore \\
{A^{4}} \ar[rr]
\ar[ur]^-{1{\times}1{\times}m}
\save \POS?="11m" \restore
\ar[dr]_-{m{\times}1{\times}1}
\save \POS?="m11" \restore
&& {A^{3}}
\ar[ur] \save \POS?="1m" \restore
\ar[dr] \save \POS?="m1" \restore
&& {A} \\
& {A^{3}} \ar[rr]_-{m{\times}1}
&& {A^{2}} \ar[ur]_-{m}
\POS "1m"; "11m" **@{}
? /-.5em/ \ar@{=>}_{1{\times}{\alpha}} ? /.5em/
\POS "m11"; "m1" **@{}
? /-.5em/ \ar@{=>}^{{\alpha}{\times}1} ? /.5em/
\POS "m1"; "m" **@{}
? /-.5em/ \ar@{=>}^{\alpha} ? /.5em/
}}}="lhs";
(1,0)*-!L{\xybox{
\xymatrix @C=1em {
& {A^{3}} \ar[rr]^-{1{\times}m}
\ar[dr] \save \POS?="m1" \restore
&& {A^{2}} \ar[dr]^-{m}
\save \POS?="mtop" \restore \\
{A^{4}}
\ar[ur]^-{1{\times}1{\times}m}
\ar[dr]_-{m{\times}1{\times}1}
\save \POS?="m11" \restore
&& {A^{2}} \ar[rr]
&& {A} \\
& {A^{3}} \ar[rr]_-{m{\times}1}
\ar[ur] \save \POS?="1m" \restore
&& {A^{2}} \ar[ur]_-{m}
\save \POS?="mbot" \restore
\POS "m11"; "m1" **@{}
? /-.35em/ \ar@{=} ? /.35em/
\POS "m1"; "mtop" **@{}
? /-.5em/ \ar@{=>}^{\alpha} ? /.5em/
\POS "mbot"; "1m" **@{}
? /-.5em/ \ar@{=>}_{\alpha} ? /.5em/
}}}="rhs"
**@{} ? /-.5em/ \ar@{=} ? /.5em/
\end{xy} \]
\[ \begin{xy}
(0,0);<2em,0em>:<0em,2em>::
(-1,0)*-!R{\xybox{
\xymatrix @C=1em {
& {A^{2}} \ar[rr]^-{1}
&& {A^{2}} \ar[dr]^-{m}
\save \POS?="m" \restore \\
{A^{2}} \ar[rr]^-{1}
\ar[ur]^-{1}
\save \POS?="11m" \restore
\ar[dr]_-{1{\times}i{\times}1}
\save \POS?="m11" \restore
&& {A^{2}}
\ar[ur]|{1} \save \POS?="1m" \restore
\ar[dr]|{1} \save \POS?="m1" \restore
&& {A} \\
& {A^{3}} \ar[rr]_-{m{\times}1}
&& {A^{2}} \ar[ur]_-{m}
\POS "1m"; "11m" **@{}
? /-.35em/ \ar@{=} ? /.35em/
\POS "m11"; "m1" **@{}
? /-.5em/ \ar@{=>}^{{\rho}{\times}1} ? /.5em/
\POS "m1"; "m" **@{}
? /-.35em/ \ar@{=} ? /.35em/
}}}="lhs";
(1,0)*-!L{\xybox{
\xymatrix @C=1em {
& {A^{2}} \ar[rr]^-{1}
\ar[dr]|{1} \save \POS?="m1" \restore
&& {A^{2}} \ar[dr]^-{m}
\save \POS?="mtop" \restore \\
{A^{2}}
\ar[ur]^-{1}
\ar[dr]_-{1{\times}i{\times}1}
\save \POS?="m11" \restore
&& {A^{2}} \ar[rr]^-{m}
&& {A} \\
& {A^{3}} \ar[rr]_-{m{\times}1}
\ar[ur]|{1{\times}m} \save \POS?="1m" \restore
&& {A^{2}} \ar[ur]_-{m}
\save \POS?="mbot" \restore
\POS "m11"; "m1" **@{}
? /-.5em/ \ar@{=>}^{1{\times}{\lambda}} ? /.5em/
\POS "m1"; "mtop" **@{}
? /-.35em/ \ar@{=} ? /.35em/
\POS "mbot"; "1m" **@{}
? /-.5em/ \ar@{=>}_{\alpha} ? /.5em/
}}}="rhs"
**@{} ? /-.5em/ \ar@{=} ? /.5em/
\end{xy} \]
A \emph{monoid} in $\ca K$ is a pseudo-monoid
for which the two-cells in the above definition
are identities.
\end{definition}
\begin{definition}
Let $(A,i,m,\alpha,\lambda,\rho))$ and $(A',i',m',{\alpha}',{\lambda}',{\rho}'))$
be pseudo monoids.
A \emph{strong monoidal structure} for a $1$-cell
$f:A{\rightarrow}A'$ consists of invertible $2$-cells
\[ \xymatrix{1 \ar[r]^-{i} \ar[d]_{1} 
\save \POS?="codph0" \restore
& A \ar[d]|{f} 
\save \POS?="domph0" \restore 
\save \POS?="domph2" \restore
& {A{\times}A} \ar[l]_-{m} \ar[d]^{f{\times}f} 
\save \POS?="codph2" \restore \\
1 \ar[r]_-{i'} & A' & {A'{\times}A'} \ar[l]^-{m'}
\save \POS"domph0";"codph0"**@{}
?  /-0.5em/ \ar@{=>}_{{\phi}_{0}} ?  /0.5em/ \restore
\save \POS"domph2";"codph2"**@{}
?  /-0.5em/ \ar@{=>}^{{\phi}_{2}} ?  /0.5em/ \restore
} \]
in $\ca K$ satisfying
the following three axioms:
\[ \begin{xy}
(0,0);<2em,0em>:<0em,2em>::
(-1,0)*-!R{\xybox{
\xymatrix @C=1em {
& {{A'}^{3}} \ar[rr]^-{m'{\times}1}
&& {{A'}^{2}} \ar[dr]^-{m'}
\save \POS?="m" \restore \\
{A^{3}} \ar[rr]
\ar[ur]^-{f}
\save \POS?="11m" \restore
\ar[dr]_-{1{\times}m}
\save \POS?="m11" \restore
&& {A^{2}}
\ar[ur] \save \POS?="1m" \restore
\ar[dr] \save \POS?="m1" \restore
&& {A'} \\
& {A^{2}} \ar[rr]_-{m}
&& {A} \ar[ur]_-{f}
\POS "1m"; "11m" **@{}
? /-.5em/ \ar@{=>}_{{\phi}_{2}{\times}f} ? /.5em/
\POS "m11"; "m1" **@{}
? /-.5em/ \ar@{=>}^{\alpha} ? /.5em/
\POS "m1"; "m" **@{}
? /-.5em/ \ar@{=>}^{{\phi}_{2}} ? /.5em/
}}}="lhs";
(1,0)*-!L{\xybox{
\xymatrix @C=1em {
& {{A'}^{3}} \ar[rr]^-{m'{\times}1}
\ar[dr] \save \POS?="m1" \restore
&& {A'^{2}} \ar[dr]^-{m'}
\save \POS?="mtop" \restore \\
{A^{3}}
\ar[ur]^-{f^{3}}
\ar[dr]_-{1{\times}m}
\save \POS?="m11" \restore
&& {A'^{2}} \ar[rr]
&& {A'} \\
& {A^{2}} \ar[rr]_-{m}
\ar[ur] \save \POS?="1m" \restore
&& {A} \ar[ur]_-{f}
\save \POS?="mbot" \restore
\POS "m11"; "m1" **@{}
? /-.5em/ \ar@{=>}^{f{\times}{\phi}_{2}} ? /.5em/
\POS "m1"; "mtop" **@{}
? /-.5em/ \ar@{=>}^{a'} ? /.5em/
\POS "mbot"; "1m" **@{}
? /-.5em/ \ar@{=>}_{{\phi}_{2}} ? /.5em/
}}}="rhs"
**@{} ? /-.5em/ \ar@{=} ? /.5em/
\end{xy} \]    
\[ \begin{xy}
(0,0);<2em,0em>:<0em,2em>::
(-1,0)*-!R{\xybox{
\xymatrix @C=1em {
& {A'} \ar[rr]^-{1}
&& {A'} \ar[dr]^-{1}
\save \POS?="m" \restore \\
{A} \ar[rr]^-{1}
\ar[ur]^-{f}
\save \POS?="11m" \restore
\ar[dr]_-{1{\times}i}
\save \POS?="m11" \restore
&& {A}
\ar[ur]|{f} \save \POS?="1m" \restore
\ar[dr]|{1} \save \POS?="m1" \restore
&& {A'} \\
& {A^{2}} \ar[rr]_-{m}
&& {A} \ar[ur]_-{f}
\POS "1m"; "11m" **@{}
? /-.35em/ \ar@{=} ? /.35em/
\POS "m11"; "m1" **@{}
? /-.5em/ \ar@{=>}^{\rho} ? /.5em/
\POS "m1"; "m" **@{}
? /-.35em/ \ar@{=} ? /.35em/
}}}="lhs";
(1,0)*-!L{\xybox{
\xymatrix @C=1em {
& {A'} \ar[rr]^-{1}
\ar[dr]|{1{\times}i'} \save \POS?="m1" \restore
&& {A'} \ar[dr]^-{1}
\save \POS?="mtop" \restore \\
{A}
\ar[ur]^-{f}
\ar[dr]_-{1{\times}i}
\save \POS?="m11" \restore
&& {A'^{2}} \ar[rr]^-{m'}
&& {A'} \\
& {A^{2}} \ar[rr]_-{m}
\ar[ur]|{f^{2}} \save \POS?="1m" \restore
&& {A} \ar[ur]_-{f}
\save \POS?="mbot" \restore
\POS "m11"; "m1" **@{}
? /-.5em/ \ar@{=>}^{f{\times}{\phi}_{0}} ? /.5em/
\POS "m1"; "mtop" **@{}
? /-.5em/ \ar@{=>}^{{\rho}'} ? /.5em/
\POS "mbot"; "1m" **@{}
? /-.5em/ \ar@{=>}_{{\phi}_{2}} ? /.5em/
}}}="rhs"
**@{} ? /-.5em/ \ar@{=} ? /.5em/
\end{xy} \]
\[ \begin{xy}
(0,0);<2em,0em>:<0em,2em>::
(-1,0)*-!R{\xybox{
\xymatrix @C=1em {
& {A'} \ar[rr]^-{1}
&& {A'} \ar[dr]^-{1}
\save \POS?="m" \restore \\
{A} \ar[rr]^-{1}
\ar[ur]^-{f}
\save \POS?="11m" \restore
\ar[dr]_-{i{\times}1}
\save \POS?="m11" \restore
&& {A}
\ar[ur]|{f} \save \POS?="1m" \restore
\ar[dr]|{1} \save \POS?="m1" \restore
&& {A'} \\
& {A^{2}} \ar[rr]_-{m}
&& {A} \ar[ur]_-{f}
\POS "1m"; "11m" **@{}
? /-.35em/ \ar@{=} ? /.35em/
\POS "m11"; "m1" **@{}
? /-.5em/ \ar@{=>}^{\lambda} ? /.5em/
\POS "m1"; "m" **@{}
? /-.35em/ \ar@{=} ? /.35em/
}}}="lhs";
(1,0)*-!L{\xybox{
\xymatrix @C=1em {
& {A'} \ar[rr]^-{1}
\ar[dr]|{i'{\times}1} \save \POS?="m1" \restore
&& {A'} \ar[dr]^-{1}
\save \POS?="mtop" \restore \\
{A}
\ar[ur]^-{f}
\ar[dr]_-{1{\times}i}
\save \POS?="m11" \restore
&& {A'^{2}} \ar[rr]^-{m'}
&& {A'} \\
& {A^{2}} \ar[rr]_-{m}
\ar[ur]|{f^{2}} \save \POS?="1m" \restore
&& {A} \ar[ur]_-{f}
\save \POS?="mbot" \restore
\POS "m11"; "m1" **@{}
? /-.5em/ \ar@{=>}^{{\phi}_{0}{\times}f} ? /.5em/
\POS "m1"; "mtop" **@{}
? /-.5em/ \ar@{=>}^{{\lambda}'} ? /.5em/
\POS "mbot"; "1m" **@{}
? /-.5em/ \ar@{=>}_{{\phi}_{2}} ? /.5em/
}}}="rhs"
**@{} ? /-.5em/ \ar@{=} ? /.5em/
\end{xy} \]
The strong monoidal morphism $(f,\phi_0,\phi_2)$
is said to be strict, when $\phi_0$ and $\phi_2$ are identities.
\end{definition}
\begin{definition}
Let $(f,\phi_0,\phi_2)$ and $(f',{\phi_0}',{\phi_2}')$
be strong monoidal morphisms.
A $2$-cell $\psi:f{\implies}f'$ is a \emph{monoidal $2$-cell} when
\[ \begin{xy}
(0,0);<2em,0em>:<0em,2em>::
(-1,0)*-!R{\xybox{
*{{\phi}_{0}}
}}="lhs";
(1,0)*-!L{\xybox{
\xymatrix{1 \ar[r]^-{i} \ar@/_1pc/[d]_{1}
\save \POS?="codphi0" \restore
& A \ar@/_1pc/[d]_{f'}
\save \POS?="domphi0" \restore
\save \POS?="codpsi" \restore
\ar@/^1pc/[d]^{f}
\save \POS?="dompsi" \restore \\
1 \ar[r]_{i'} & A'
\POS "domphi0"; "codphi0" **@{}
? /-.5em/ \ar@{=>}_{\phi_{0}} ? /.5em/
\POS "dompsi"; "codpsi" **@{}
? /-.5em/ \ar@{=>}_{\psi} ? /.5em/ }
}}="rhs"
**@{} ? /-.5em/ \ar@{=} ? /.5em/
\end{xy} \]
\[ \begin{xy}
(0,0);<2em,0em>:<0em,2em>::
(-1,0)*-!R{\xybox{
\xymatrix{A \ar@/_1pc/[d]_{f}
\save \POS?="dompsi" \restore
\ar@/^1pc/[d]^{f'}
\save \POS?="codpsi" \restore
\save \POS?="domphip2" \restore
& {A^{2}} \ar[l]_-{m}
\ar@/^1pc/[d]^{f'^{2}}
\save \POS?="codphip2" \restore \\
{A'} & {A'^{2}} \ar[l]^-{m'}
\POS "dompsi"; "codpsi" **@{}
? /-.5em/ \ar@{=>}^{\psi} ? /.5em/
\POS "domphip2"; "codphip2" **@{}
? /-.5em/ \ar@{=>}^{{\phi}'_{2}} ? /.5em/
}
}}="lhs";
(1,0)*-!L{\xybox{
\xymatrix{A \ar@/_1pc/[d]_{f}
\save \POS?="domphi2" \restore
& {A^{2}} \ar[l]_-{m}
\ar@/_1pc/[d]_{f^{2}}
\save \POS?="dompsi2" \restore
\save \POS?="codphi2" \restore
\ar@/^1pc/[d]^{f'^{2}}
\save \POS?="codpsi2" \restore \\
{A'} & {A'^{2}} \ar[l]^-{m'}
\POS "dompsi2"; "codpsi2" **@{}
? /-.5em/ \ar@{=>}^{{\psi}^{2}} ? /.5em/
\POS "domphi2"; "codphi2" **@{}
? /-.5em/ \ar@{=>}^{\phi_{2}} ? /.5em/
}
}}="rhs"
**@{} ? /-.5em/ \ar@{=} ? /.5em/
\end{xy} \]
\end{definition}
With the evident compositions,
one defines the $2$-category $\PsMon {\ca K}$
to consist of pseudo monoids,
strong monoidal morphisms and monoidal $2$-cells.
The locally full sub-$2$-category of $\PsMon {\ca K}$
consisting of strict monoids and strict monoid morphisms
is denoted as $\Mon {\ca K}$.
\begin{example}
The $2$-category $\PsMon {\CAT}$ consists of
monoidal categories, strong monoidal functors, and monoidal natural
transformations in the usual sense.
\end{example}
Given a pseudo monoid structure $(i,m,\alpha,\lambda,\rho))$ on $A \in \ca K$,
then by composition, for each $X \in \ca K$, the hom category $\ca K(X,A)$
obtains a monoidal category structure,
and these monoidal structures are $2$-natural in $X$.
This is true since the monoidal structure of $\ca K$ is cartesian product,
and representable $2$-functors preserve products.
On the other hand, by the $\CAT$-enriched yoneda lemma,
monoidal category structures on the homs $\ca K(X,A)$,
$2$-natural in $X$, determines a pseudo monoid structure on $A$.
This definition of a pseudo monoid structure on $A$ via the homs $\ca K(X,A)$
is called the \emph{representable} definition.
Strong monoidal morphisms and monoidal $2$-cells
can be defined representably in the same way.
Note that the forgetful $2$-functor
\[ \xymatrix{{\PsMon {\ca K}} \ar[r] & {\ca K}} \]
can easily be seen to create products.
\begin{examples}
\begin{enumerate}
\item  Let $\ca E$ be a category with finite products,
and let $\ca K$ be $\ca E$ regarded
as a locally discrete $2$-category.
Then a pseudo monoid is just a monoid in $\ca E$ in the usual sense.
\item  Let $\ca K$ be $\PsMon {\CAT}$.
By \cite{JS93}, a pseudo monoid in $\ca K$
is a braided monoidal category.
\label{Psm-Mon}
\item  Let $\ca K$ be $\PsMon {\PsMon {\CAT}}$.
By \cite{JS93}, a pseudo monoid in $\ca K$
is a symmetric monoidal category.
\item  Let $\ca K$ be $\PsMon {\PsMon {\PsMon {\CAT}}}$.
By \cite{JS93}, a pseudo monoid in $\ca K$
is also a symmetric monoidal category.
\end{enumerate}
\end{examples}
In fact, the results of \cite{JS93} alluded
to in the above examples actually assert $2$-equivalences,
that is, equivalences in the $\CAT$-enriched sense,
between the appropriate $2$-categories.
In particular, the forgetful $2$-functor
\[ \xymatrix{{\PsMon {\PsMon {\PsMon {\CAT}}}} \ar[r]
& {\PsMon {\PsMon {\CAT}}}} \]
is a $2$-equivalence.
Since the $2$-categories $\PsMon {\ca K}$ can be defined representably,
one immediately obtains the following ``Eckmann Hilton'' stabilisation result.
\begin{proposition}\label{EH2}
Let $\ca K$ be a $2$-category with finite products.
Then the forgetful $2$-functor
\[ \xymatrix{{\PsMon {\PsMon {\PsMon {\ca K}}}} \ar[r]
& {\PsMon {\PsMon {\ca K}}}} \]
is a $2$-equivalence.
\end{proposition}
%

\section{Monoidal pseudo algebras}\label{mon-psalg}

For this section, let $T$ be a $2$-monad on a $2$-category $\ca K$
with finite products.
It is easily seen that both the forgetful $2$-functors
\[ \begin{array}{rccl}
{\PsAlg T}{\rightarrow}{\ca K}
&&& {\NormPsAlg T}{\rightarrow}{\ca K}
\end{array}
\]
create products, and so in particular, $\NormPsAlg T$
has finite products.
\begin{definition}
A \emph{monoidal pseudo $T$-algebra}
is a pseudo monoid in $\NormPsAlg T$.
\end{definition}
Unpacking this definition,
one finds that a monoidal pseudo $T$-algebra consists of
\begin{itemize}
\item  an object $A \in \ca K$.
\item  a normal pseudo $T$-algebra structure $(a,\alpha)$
on $A$.
\item  a pseudo monoid structure $(i,m,\beta,\lambda,\rho)$
on $A$.
\item  an invertible $2$-cell $\overline{i}$
which provides $i$ with a strong $T$-morphism structure.
\item  an invertible $2$-cell $\overline{m}$
which provides $m$ with a strong $T$-morphism structure.
\item  the $2$-cells $\beta$, $\lambda$ and $\rho$
satisfy the $T$-algebra $2$-cell axiom.
\end{itemize}
We shall refer to this monoidal pseudo algebra
by the ordered $8$-tuple $(A,a,i,m,\alpha,\beta,\lambda,\rho)$.
\begin{examples}\label{ex:mon-psalg}
\begin{enumerate}
\item  For $T$ as in (\ref{Monad-ex})(\ref{loc-discrete}),
a monoidal pseudo algebra is a monoid in $\Alg T$.
\item  For $T = 1_{\CAT}$
a monoidal pseudo algebra is a monoidal category.
More generally, for $T = 1_{\ca K}$ as in (\ref{Monad-ex})(\ref{id}),
a monoidal pseudo algebra is a pseudo monoid.
\item  For $T = \ca M$ as in (\ref{Monad-ex})(\ref{Mon-on-Cat}),
a monoidal pseudo algebra is a braided monoidal category.
Abstractly, this follows from (\ref{EH2}), because $\NormPsAlg {\ca M}$
is $2$-equivalent to the $2$-category of monoidal categories,
strong monoidal functors, and monoidal natural transformations.
In this formalism the braiding arises from $\overline{m}$.
In more detail, denote the object map of $m$ by $m(x,y) = x \tensor_0 y$,
and the object map of $a$ by
$a(x_0,...,x_{n-1}) = x_0 \tensor_1 ... \tensor_1 x_{n-1}$.
Then $\overline{m}$ is an invertible $2$-cell
\[ \xymatrix{{{\ca M}(A{\times}A)} \ar[r] \ar[d]_{M(m)}
\save \POS?="dom" \restore & {{\ca M}(A){\times}{\ca M}(A)}
\ar[r]^-{a{\times}a} & {A{\times}A} \ar[d]^{m}
\save \POS?="cod" \restore \\
{{\ca M}(A)} \ar[rr]_-{a} && A
\POS "dom"; "cod" **@{};
? /-.5em/ \ar@{=>}^{\overline{m}} ? /.5em/ } \]
where ${\ca M}(A{\times}A){\rightarrow}{\ca M}(A){\times}{\ca M}(A)$
is the canonical comparison.
So the component of $\overline{m}$ at $((x_0,y_0),...,(x_{n-1},y_{n-1}))$
is an isomorphism
\[ \xymatrix{(x_0{\tensor_0}y_0){\tensor_1}...{\tensor_1}(x_{n-1}{\tensor_0}y_{n-1})
\ar[d] \\
(x_0{\tensor_1}...{\tensor_1}x_{n-1}){\tensor_0}
(y_0{\tensor_1}...{\tensor_1}y_{n-1})} \]
which in the current context deserves to be called the braiding.
Writing $I_0$ for the unit for $\tensor_0$,
the components of $\overline{i}$ are isomorphisms
\[ {\underbrace{I_0{\tensor_1}...{\tensor_1}I_0}_n} {\rightarrow} {I_0} \]
which in the case $n=0$, gives an isomorphism $I_1{\iso}I_0$,
where $I_1$ is the unit for $\tensor_1$.
Furthermore $x{\tensor_1}y{\iso}x{\tensor_0}y$ is obtained as:
\[ \begin{array}{rcccl}
x{\tensor_1}y & \iso & (x{\tensor_0}I_0)\tensor_1(I_0{\tensor_0}y)
& \iso & (x{\tensor_1}I_0)\tensor_0(I_0{\tensor_1}y) \\
& \iso & (x{\tensor_1}I_1)\tensor_0(I_1{\tensor_1}y)
& \iso & x{\tensor_0}y
\end{array}, \]
and a braiding in the usual sense is obtained as:
\[ \begin{array}{rcccl}
x{\tensor_0}y & \iso & x{\tensor_1}y
& \iso & (I_0{\tensor_0}x)\tensor_1(y{\tensor_0}I_0) \\
& \iso & (I_0{\tensor_1}y)\tensor_0(x{\tensor_1}I_0)
& \iso & (I_1{\tensor_1}y)\tensor_0(x{\tensor_1}I_1) \\
& \iso & y{\tensor_0}x
\end{array}. \]
These isomorphisms encode the
Eckmann-Hilton argument (see \cite{Mac71}, pg 45, exercise 5).
\label{monMalg}
\item  For $T = \ca B$ as in (\ref{Monad-ex})(\ref{BrMon-on-Cat}),
a monoidal pseudo algebra is a braided monoidal category.
Abstractly, this follows from (\ref{EH2}), because $\NormPsAlg {\ca B}$
is $2$-equivalent to the $2$-category of braided monoidal categories,
braided strong monoidal functors, and monoidal natural transformations.
The more explicit analysis here only differs from the previous example
in that the action $a$ already carries the information of a braiding
for $\tensor_1$.
The naturality of $\overline{m}$ ensures that
the braiding encoded by it coincides with that described by $a$,
and forces it to be a symmetry.
\item  Similarly for $T = \ca S$ as in (\ref{Monad-ex})(\ref{BrMon-on-Cat}),
a monoidal pseudo algebra is a symmetric monoidal category.
The more explicit analysis here only differs from the previous example
in that the action $a$ already carries the information of a symmetry,
and so $\overline{m}$ encodes no new information.
\label{symmon}
\end{enumerate}
\end{examples}
Further examples relevant to higher dimensional algebra
will be considered in section(\ref{higher-sym-ops}).
We shall now express the pseudo monoid part of a monoidal pseudo
algebra representably, to facilitate the general operad definition.
To this end, let $(A,a,i,m,\alpha,\beta,\lambda,\rho)$
be a monoidal pseudo $T$-algebra.
First, we note that the unit object
$i:X{\rightarrow}A$ of the monoidal category ${\ca K}(X,A)$ is the composite
\[ \xymatrix{X \ar[r]^-{!} & 1 \ar[r]^-{i} & A} \]
in $\ca K$, where $!$ here denotes the unique map into the terminal object.
Moreover, given objects $x$ and $y$ of ${\ca K}(X,A)$,
their tensor product $x{\tensor}y$ is the composite
\[ \xymatrix{X \ar[r]^-{(x,y)} & {A^2} \ar[r]^-{m} & A} \]
in $\ca K$. Note that if $z:Z{\rightarrow}X$,
then $(x{\tensor}y)z = xz{\tensor}yz$ by the naturality of $\tensor$.
Similarly one can express the rest of the pseudo monoid data $(\beta,\lambda,\rho)$
representably.
One can write $\overline{i}:aT(i){\rightarrow}i$
for the $2$-cell
\[ \xymatrix{{T1} \ar[r]^-{!} \ar[d]_{T(i)}
\save \POS?="domibar" \restore & 1 \ar[d]^{i}
\save \POS?="codibar" \restore \\ T(A) \ar[r]_-{a} & A
\POS "domibar"; "codibar" **@{};
? /-.5em/ \ar@{=>}^{\overline{i}} ? /.5em/ } \]
which provides $i$'s strong $T$-morphism structure.
As for $\overline{m}$,
given objects $x$ and $y$ of ${\ca K}(X,A)$,
we shall write
\[ \xymatrix{{aT(x{\tensor}y)} \ar[r]^-{\overline{m}_{x,y}}
& {aT(x){\tensor}aT(y)}} \]
for the composite
\[ \xymatrix{{T(X)} \ar[r]^-{T(x,y)} & {T(A{\times}A)} \ar[r]^-{\pi} \ar[d]_{T(m)}
\save \POS?="dommbar" \restore
& {T(A){\times}T(A)} \ar[r]^-{a{\times}a} & {A{\times}A} \ar[d]^{m}
\save \POS?="codmbar" \restore \\ & T(A) \ar[rr]_-{a} && A
\POS "dommbar"; "codmbar" **@{};
? /-.5em/ \ar@{=>}^{\overline{m}} ? /.5em/ } \]
where $\pi$ is the canonical comparison,
and $\overline{m}$ is $m$'s strong $T$-morphism structure.
When the context is clear we shall drop the subscripts and write
\[ \overline{m}:aT(x{\tensor}y){\rightarrow}aT(x){\tensor}aT(y). \]
In light of this notation,
the strong $T$-morphism axioms for $\overline{m}$,
and the $T$-algebra $2$-cell axioms for $\beta$, $\lambda$
and $\rho$, can be restated as follows.
\begin{proposition}\label{rep-axioms}
\begin{enumerate}
\item  $\forall x,y \in {\ca K}(X,A)$,
\[ \xymatrix{{aT(a)T^2(x{\tensor}y)} \ar[rr]^-{{\alpha}T^2(x{\tensor}y)}
\ar[d]_{aT(\overline{m})} && {aT(x{\tensor}y)\mu_{TA}} 
\ar[dd]^{\overline{m}\mu_{TA}} \\
{aT(aT(x){\tensor}aT(y))} \ar[d]_{\overline{m}} \\
{aT(a)T^2(x){\tensor}aT(a)T^2(y)} \ar[rr]_-{{\alpha}T^2(x){\tensor}{\alpha}T^2(y)}
&& {aT(x)\mu_{TA}{\tensor}aT(y)\mu_{TA}}} \]
commutes in ${\ca K}(X,A)$.
\label{mbar-Tmor}
\item  $\forall x,y \in {\ca K}(X,A)$, $\overline{m}_{x,y}\eta_X = 1_{x{\tensor}y}$.
\item  $\forall x,y,z \in {\ca K}(X,A)$,
\[ \xymatrix{{aT(x{\tensor}(y{\tensor}z))}
\ar[rr]^-{aT(\beta)} \ar[d]_{\overline{m}} && {aT((x{\tensor}y){\tensor}z)}
\ar[d]^{\overline{m}} \\ {aT(x){\tensor}aT(y{\tensor}z)}
\ar[d]_{\id\tensor\overline{m}} && {aT(x{\tensor}y){\tensor}aT(z)}
\ar[d]^{\overline{m}} \\ {aT(x){\tensor}(aT(y){\tensor}aT(z))}
\ar[rr]_-{\beta} && {(aT(x){\tensor}aT(y)){\tensor}aT(z)}} \]
commutes in ${\ca K}(X,A)$.
\item  $\forall x \in {\ca K}(X,A)$,
\[ \begin{xy}
(0,0);<2em,0em>:<0em,2em>::
(-1,0)*+!R{\xybox{
\xymatrix{{aT(i{\tensor}x)} \ar[d]_{aT(\lambda)} \ar[r]^-{\overline{m}}
& {aT(i){\tensor}aT(x)} \ar[d]^{\overline{i}\tensor\id} \\
aT(x) & {i{\tensor}aT(x)} \ar[l]^-{\lambda}}}},
(1,0)*+!L{\xybox{
\xymatrix{{aT(x{\tensor}i)} \ar[d]_{aT(\rho)} \ar[r]^-{\overline{m}}
& {aT(x){\tensor}aT(i)} \ar[d]^{\id\tensor\overline{i}} \\
aT(x) & {aT(x){\tensor}i} \ar[l]^-{\rho}}}}
\end{xy} \]
commute in $\ca K(X,A)$.
\end{enumerate}
\end{proposition}
When writing diagrams such as those in (\ref{rep-axioms}),
notice that there are situations when objects can be expressed
in more than one way. For instance in (\ref{rep-axioms})(\ref{mbar-Tmor})
we have $a\mu_AT^2(x{\tensor}y) = aT(x{\tensor}y)\mu_{TA}$
by the naturality of $\mu$,
although in that diagram we have only recorded $aT(x{\tensor}y)\mu_{TA}$.
In similar situations below,
we shall just choose one description of a given object without further comment
when there is little risk of confusion.

\section{Operads}\label{operads}

With the language of monoidal pseudo algebras
at our disposal,
we are now able to present our general operad definition.
\begin{definition}
Let $(T,\eta,\mu)$ be a $2$-monad on a $2$-category $\ca K$ with finite products,
and let $(A,a,i,m,\alpha,\beta,\lambda,\rho)$
be a monoidal pseudo $T$-algebra.
A \emph{$T$-operad} $(p,\iota,\sigma)$ in $A$
consists of a $1$-cell $p:T(1){\rightarrow}A$,
together with $2$-cells $\iota$ and $\sigma$
\[ \begin{xy}
(0,0);<2em,0em>:<0em,2em>::
(-1,0)*+!R{\xybox{
\xymatrix{
1 \ar[rr]^-{\eta_1} \ar[dr]_{i} \save \POS?="domiot" \restore
&& {T1} \ar[dl]^{p} \save \POS?="codiot" \restore \\ & A
\POS "domiot";"codiot" **@{}
? /-.5em/ \ar@{=>}^{\iota} ? /.5em/ }}},
(1,0)*+!L{\xybox{
\xymatrix{
{T^2{1}} \ar[dr]_{pT(!){\tensor}aT(p)} \save \POS?="domsig" \restore
\ar[rr]^-{\mu_1}
&& {T1} \ar[dl]^{p}  \save \POS?="codsig" \restore \\
& A
\POS "domsig";"codsig" **@{}
? /-.5em/ \ar@{=>}^{\sigma} ? /.5em/ }}}
\end{xy} \]
such that
\[ \begin{xy}
(0,0);<2em,0em>:<0em,2em>::
(-1,0)*+!R{\xybox{
\xymatrix{{i{\tensor}p}
\ar[r]^-{{\iota}!{\tensor}\id}
\ar[dr]_{\lambda}
& {pT(!)\eta_{T1}{\tensor}p}
\ar[d]^{\sigma\eta_{T1}} \\ & p}}},
(1,0)*+!L{\xybox{
\xymatrix{{p{\tensor}aT(i)} \ar[d]_{\id{\tensor}\overline{i}}
\ar[r]^-{\id\tensor{a}T(\iota)}
& {p{\tensor}aT(p)T(\eta_1)}
\ar[d]^{{\sigma}T(\eta_1)} \\
{p{\tensor}i} \ar[r]_-{\rho} & p}}}
\end{xy} \]
commute in ${\ca K}(T1,A)$ and
\[ \xymatrix{{pT(!){\tensor}aT(pT(!){\tensor}aT(p))}
\ar[r]^-{\id\tensor\overline{m}}
\ar[d]_{\id{\tensor}aT(\sigma)}
& {pT(!)\tensor(aT(p)T^2(!){\tensor}aT(a)T^2(p))}
\ar[d]^{\beta} \\
{pT(!){\tensor}aT(p)T(\mu_1)}
\ar[d]_{{\sigma}T(\mu_1)}
& {(pT(!){\tensor}aT(p)T^2(!)){\tensor}aT(a)T^2(p)}
\ar[d]^{{\sigma}T^2(!){\tensor}{\alpha}T^2(p)} \\
{p{\mu_1}T(\mu_1)}
& {pT(!)\mu_{T1}{\tensor}aT(p)\mu_{T1}}
\ar[l]^-{\sigma\mu_{T1}}} \]
commutes in ${\ca K}(T^31,A)$.
\end{definition}
\begin{definition}
Let $(T,\eta,\mu)$ be a $2$-monad on a $2$-category $\ca K$ with finite products,
and let $(A,a,i,m,\alpha,\beta,\lambda,\rho)$
be a monoidal pseudo $T$-algebra.
A \emph{morphism} $(p,\iota,\sigma){\rightarrow}(p',{\iota}',{\sigma}')$
of $T$-operads in $A$ consists of a $2$-cell $\phi:p{\implies}p'$
such that
\[ \begin{xy}
(0,0);<2em,0em>:<0em,2em>::
(-1,0)*+!R{\xybox{
\xymatrix{i \ar[r]^-{\iota} \ar[dr]_{{\iota}'} & {p\eta_1}
\ar[d]^{\phi\eta_1} \\ & {p'\eta_1}}}},
(1,0)*+!L{\xybox{
\xymatrix{{pT(!){\tensor}aT(p)} \ar[r]^-{\sigma}
\ar[d]_{{\phi}T(!){\tensor}aT(\phi)} & {p\mu_1} \ar[d]^{\phi\mu_1} \\
{p'T(!){\tensor}aT(p')} \ar[r]_-{{\sigma}'} & {p'\mu_1}}}}
\end{xy} \]
commute ${\ca K}(1,A)$ and ${\ca K}(T^21,A)$ respectively.
\end{definition}
With the evident composition, one obtains the category
$\Op T A$ of $T$-operads in the monoidal pseudo $T$-algebra $A$,
and a forgetful functor ${\Op T A}{\rightarrow}{\ca K}(T1,A)$.
\begin{definition}
Let $(p,\iota,\sigma)$ be a $T$-operad in $A$.
A $p$-algebra $(x,\overline{x})$
consists of a one-cell $x:1{\rightarrow}A$ and a $2$-cell
\[ \xymatrix{{T1} \ar[rr]^-{!} \ar[dr]_{p{\tensor}aT(x)}
\save \POS?="domxbar" \restore
&& 1 \ar[dl]^{x} \save \POS?="codxbar" \restore \\ & A
\POS "domxbar";"codxbar" **@{}
? /-.5em/ \ar@{=>}^{\overline{x}} ? /.5em/ } \]
such that
\[ \xymatrix{{i{\tensor}x} \ar[rr]^-{\iota\tensor\id} \ar[dr]_{\lambda}
&& {p{\eta_1}{\tensor}x} \ar[dl]^{\overline{x}\eta_1} \\ & x} \]
commutes in ${\ca K}(1,A)$ and
\[ \xymatrix{{pT(!){\tensor}aT(p{\tensor}aT(x))}
\ar[r]^-{\id\tensor\overline{m}}
\ar[d]_{\id{\tensor}aT(\overline{x})}
& {pT(!)\tensor(aT(p){\tensor}aT(a)T^2(x))}
\ar[d]^{\beta} \\
{pT(!){\tensor}aT(x)T(!)}
\ar[d]_{\overline{x}T(!)}
& {(pT(!){\tensor}aT(p)){\tensor}aT(a)T^2(x)}
\ar[d]^{{\sigma}{\tensor}{\alpha}T^2(x)} \\
{x!}
& {p\mu_1{\tensor}aT(x)\mu_1}
\ar[l]^-{\overline{x}\mu_1}} \]
commutes in ${\ca K}(T^21,A)$.
\end{definition}
\begin{definition}
Let $(p,\iota,\sigma)$ be a $T$-operad in $A$.
A $p$-algebra morphism
\[ f:(x,\overline{x}){\rightarrow}(y,\overline{y}) \]
consists of $f:x{\rightarrow}y$ in ${\ca K}(1,A)$ such that
\[ \xymatrix{{p{\tensor}aT(x)} \ar[r]^-{\overline{x}} \ar[d]_{\id{\tensor}aT(f)}
& {x!} \ar[d]^{f!} \\ {p{\tensor}aT(y)} \ar[r]_-{\overline{y}} & {y!}} \]
commutes in $\ca K(T1,A)$.
\end{definition}
With the evident composition, one obtains the category
$\Alg p$ of $p$-algebras and a forgetful functor
${\Alg p}{\rightarrow}{\ca K}(1,A)$.
We shall now see how the well known operad notions are captured
by these general definitions.
\begin{examples}\label{conv-op-ex}
\begin{enumerate}
\item  For $T = 1_{\CAT}$, a $T$-operad $(p,\iota,\sigma)$ in $A$
is a monoid $M$ in the monoidal category $A$.
The underlying object of $M$ is picked out by
$p:1{\rightarrow}A$, and the unit and multiplication are provided by
$\iota$ and $\sigma$ respectively.
An object of $\Alg p$ is an object of $A$ acted on by $M$.
\item  For $T = \ca M$, a $T$-operad $(p,\iota,\sigma)$ in $A$
is a non-symmetric operad in the braided monoidal category $A$.
In more detail, $p:{\ca M}(1){\rightarrow}A$ is a sequence
of objects $(p_n : n \in \N)$ of $A$. The unit $\iota$
amounts to a map $i{\rightarrow}p_1$ in $A$.
As for the substitution,
there is a component of $\sigma$ for each element of ${\ca M}^2(1)$,
that is, for each finite sequence $(n_j : j \in k)$
of natural numbers. The component of $\sigma$ for this sequence
is a map
\[ \xymatrix{{p_k \tensor (p_{n_0} \tensor ... \tensor p_{n_{k-1}})}
\ar[rr] && {p_n}} \]
in $A$, where $n = \sum_{j{\in}k} n_j$.
The axioms express the usual unit and associativity
laws for substitution.
Notice how the braiding $\overline{m}$ is necessary
to express the associativity of the substitution $\sigma$.
The category $\Alg p$
is the usual category of algebras for the operad.
\item  For $T = \ca B$, a $T$-operad $(p,\iota,\sigma)$ in $A$
is a braided operad in the symmetric monoidal category $A$.
This example differs from the previous one in two respects.
The first is that the functoriality of $p:{\ca B}(1){\rightarrow}A$
amounts to equipping each $p_n$
with an action of ${\Br}_n$, the $n$-th braid group.
The second is that the naturality of $\sigma$
amounts to the substitution being equivariant
with respect to these actions.
Similarly, for a $p$-algebra $(x,\overline{x})$,
the naturality of $\overline{x}$
encodes its equivariance as an action on $x$.
\item  In the same way, for $T = \ca S$, a $T$-operad $(p,\iota,\sigma)$ in $A$
is a symmetric operad in the symmetric monoidal category $A$,
with the functorialty of $p$ encoding the symmetric group actions on the $p_n$,
and the naturality of $\sigma$ encoding the equivariance.
\end{enumerate}
\end{examples}
From the above discussion and definitions,
there are two obvious questions, important to examples, to consider:
\begin{enumerate}
\item  Under what conditions is the forgetful functor
${\Op T A}{\rightarrow}{\ca K}(T1,A)$ monadic?
Of course when this happens,
one can construct free operads.
\item  Under what conditions is the forgetful functor
${\Alg p}{\rightarrow}{\ca K}(1,A)$ monadic?
When this happens, one has a monad on $\ca K(1,A)$ associated to $p$,
whose category of Eilenberg-Moore algebras is $\Alg p$.
\end{enumerate}
In the forthcoming \cite{AnOp}, it is shown that if the monoidal pseudo
algebra in question is \emph{distributive}
in a certain sense,
then both of the above forgetful functors are monadic.
For example, a monoidal pseudo $\ca M$-algebra $\ca V$,
that is, a braided monoidal category,
is distributive in the sense of \cite{AnOp},
when it has coproducts which distribute with the tensor product of $\ca V$.

\section{Higher operads}\label{sec:higher-op}

In order to understand the motivating examples of this paper,
it is necessary to review some of the combinatorial
aspects of the globular approach to higher dimensional algebra.
For a fuller discussion, see \cite{Bat98}, \cite{WebGen},
\cite{Web}, and \cite{Lei}.
Define the category $\G$ to have natural numbers as objects,
and a generating subgraph
\[ \xymatrix{0 \ar@<1ex>[r]^{\sigma_{0}} \ar@<-1ex>[r]_{\tau_{0}} 
& 1 \ar@<1ex>[r]^{\sigma_{1}} \ar@<-1ex>[r]_{\tau_{1}}
& 2 \ar@<1ex>[r]^{\sigma_{2}} \ar@<-1ex>[r]_{\tau_{2}}
& 3 \ar@<1ex>[r]^{\sigma_{3}} \ar@<-1ex>[r]_{\tau_{3}}
& \ldots} \]
subject to the ``cosource/cotarget'' equations
${\sigma_{n+1}}{\sigma_{n}} = {\tau_{n+1}}{\sigma_{n}}$
and
${\tau_{n+1}}{\tau_{n}} = {\sigma_{n+1}}{\tau_{n}}$,
for every $n \in \N$.
The objects of the category $\Glob$, are called \emph{globular sets}.
Thus, a globular set $Z$ consists of a diagram
of sets and functions
\[ \xymatrix{Z_{0} 
& Z_{1} \ar@<-1ex>[l]_{s_{0}} \ar@<1ex>[l]^{t_{0}}
& Z_{2} \ar@<-1ex>[l]_{s_{1}} \ar@<1ex>[l]^{t_{1}}
& Z_{3} \ar@<-1ex>[l]_{s_{2}} \ar@<1ex>[l]^{t_{2}}
& \ldots \ar@<-1ex>[l]_{s_{3}} \ar@<1ex>[l]^{t_{3}}} \]
so that
${s_{n}}{s_{n+1}} = {s_{n}}{t_{n+1}}$
and
${t_{n}}{t_{n+1}} = {t_{n}}{s_{n+1}}$
for every $n \in \N$.
The elements of $Z_{n}$ are called the $n$-cells of $Z$,
and the functions $s_{n}$ and $t_{n}$
are called source and target functions.
Define $Z$ to be \emph{of dimension n}
when there are no $m$-cells for $m>n$.
All constructions that we consider below,
apply equally well to the category of $n$-globular sets,
where $\G$ is replaced by the full subcategory $\G_{(n)}$
consisting of the natural numbers $\leq n$.

Let $Z$ be a globular set.
Recall from \cite{Str91} the \emph{solid triangle order}
$\blacktriangleleft$ on the elements (of all dimensions) of Z.
Define first the relation $x \prec y$ for $x \in Z_{n}$
iff $x = s_{n}(y)$ or $t_{n-1}(x) = y$. 
Then take $\blacktriangleleft$ to be the
reflexive-transitive closure of $\prec$. 
Write $\Sol(Z)$ for the preordered set so obtained.
Observe that $\Sol$ is the object map of a functor
\[ \xymatrix@1{{\Glob} \ar[r]^-{\Sol} & {\PreOrd}} \]
where $\PreOrd$ is the category of preordered sets 
and order-preserving functions.
\begin{definition}\label{glob-card}
A \emph{globular cardinal} is a globular set $Z$
such that $\Sol(Z)$ is a non-empty finite linear order.
\end{definition}
Denote by $\Theta_0$
the full subcategory of $\Glob$
consisting of the globular cardinals.
Globular cardinals
are the pasting schemes appropriate to the Batanin
definition of weak $\omega$-category \cite{Bat98},
are analysed from the present point of view
in \cite{Web}.
In particular we have
\begin{proposition}\label{theta0-old}
\begin{enumerate}
\item  Globular cardinals are finite and connected as globular sets.
\item  All morphisms in $\Theta_0$ are monic.
\item  If $X$ is a globular cardinal, then a retraction
$X{\rightarrow}Y$ of globular sets is an isomorphism.
\end{enumerate}
\end{proposition}
Write $\Tr_n$ for the set of isomorphism classes
of globular cardinals of dimension $n$.
One of the most beautiful ideas in \cite{Bat98},
is the identification of $\Tr_n$
with $n$-stage trees, where an $n$-stage tree $T$ is
defined to be a sequence
\[ T_n {\rightarrow} ... {\rightarrow} {T_0} \]
of maps in $\Simp$, the category of finite ordinals
and monotone maps, where $T_0=1$.
Central to the Batanin approach to higher dimensional algebra
is the monad $\ca T_0$ on $\Glob$ whose algebras are strict
$\omega$-categories.
The underlying functor of this monad can be described as
\[ {\ca T_0}(X)_n = \sum_{T\in\Tr_n} \Glob(T,X), \]
and the multiplication of this monad,
which encodes the pasting of globular pasting schemes,
can be specified in terms of trees.
This monad is \emph{cartesian},
in the sense that the underlying endofunctor preserves pullbacks,
and the naturality squares for $\eta$ and $\mu$ are pullback squares.
As mentioned in the introduction,
one can then regard $\ca T_0$ as a monad on $\GLOBSET$,
and then apply $\Cat$ to obtain the cartesian $2$-monad
$\ca T$ on $\GLOBCAT$.
Normal pseudo $\ca T$-algebras can be identified
with the monoidal globular categories of \cite{Bat98}.
Their relationship is analogous to
the relation between monoidal categories defined via $k$-ary tensor products
on the one hand (normal pseudo $\ca M$-algebras),
and those defined the conventional way using binary tensor products
(pseudo monoids in $\CAT$).{\footnotemark{\footnotetext{
In fact, the corresponding $2$-equivalence
for monoidal categories,
can be seen as the
restriction of the $2$-equivalence described here,
since $1$-object, monoidal $2$-globular categories,
amount to monoidal categories described in the biased fashion,
whereas $1$-object normal pseudo $\ca T_{(2)}$-algebras
amount to unbiased monoidal categories.
Here, $\ca T_{(2)}$ is the ``truncation'' of $\ca T$ to
$[\op \G_{(2)},\CAT]$.}}}
One has $2$-functors
\[ \xymatrix{{\ca MG} \ar@<1ex>[rr]^-{F} &&
{\NormPsAlg {\ca T}} \ar@<1ex>[ll]^-{G}} \]
which can be verified directly to provide a $2$-equivalence
of $2$-categories.
Given a monoidal globular category $X$,
by making a choice of bracketting of iterated expressions,
one constructs the normal pseudo $\ca T$-algebra $F(X)$,
with the same underlying globular category.
On the other hand,
given a pseudo $\ca T$-algebra $Y$,
one obtains the monoidal globular category $G(Y)$,
with the same underlying globular category,
by considering only the nullary and binary operations,
and associated coherence data.
Using the coherence results of \cite{Bat98},
one can verify directly that $F$ and $G$ form a $2$-equivalence of $2$-categories.

Pseudo monoids in $\ca {MG}$
are particularly easy to describe:
to give $X \in [\op \G,\Cat]$ a structure of pseudo monoid
in $\ca {MG}$, is the same as giving the globular category
\[ \xymatrix{1 
& X_{0} \ar@<-1ex>[l]_{s} \ar@<1ex>[l]^{t}
& X_{1} \ar@<-1ex>[l]_{s} \ar@<1ex>[l]^{t}
& X_{2} \ar@<-1ex>[l]_{s} \ar@<1ex>[l]^{t}
& \ldots \ar@<-1ex>[l]_{s} \ar@<1ex>[l]^{t}} \]
the structure of a monoidal globular category.
Such a structure was called an augmented monoidal globular category in \cite{Bat98}.
From the discussion of the previous paragraph,
we may identify $\PsMon {\NormPsAlg {\ca T}}$
as the $2$-category of augmented monoidal globular categories.
We recall some of the main examples from \cite{Bat98}.
\begin{examples}\label{ex:mgc}
\begin{enumerate}
\item  There is a $2$-functor
\[ \Span : \CAT \rightarrow \GLOBCAT \]
for which $\Span(\ca E)_n = [\op {(\G/n)},\ca E]$.
When $\ca E$ has pullbacks, there is a canonical monoidal globular structure
on $\Span(\ca E)$, and when in addition $\ca E$ has products,
this structure is augmented, with the additional (pseudo monoid) structure
being given by pointwise cartesian product in the categories $[\op {(\G/n)},\ca E]$.
\label{span}
\item  A monoidal structure on a category $\ca V$
amounts to a monoidal $2$-globular structure on
\[ \xymatrix{1 & {\ca V} \ar@<-1ex>[l] \ar@<1ex>[l]}. \]
\label{mon}
\item  A braided monoidal structure on a category $\ca V$
amounts to a monoidal $3$-globular structure on
\[ \xymatrix{1 & 1 \ar@<-1ex>[l] \ar@<1ex>[l]
& {\ca V} \ar@<-1ex>[l] \ar@<1ex>[l]}. \]
\label{brmon}
\item  A symmetric monoidal structure on a category $\ca V$
amounts to a monoidal $(n+1)$-globular structure on
\[ \xymatrix{1 & {...} \ar@<-1ex>[l] \ar@<1ex>[l]
& 1 \ar@<-1ex>[l] \ar@<1ex>[l]
& {\ca V} \ar@<-1ex>[l] \ar@<1ex>[l]} \]
where $n{\geq}3$.
\label{symmon}
\end{enumerate}
\end{examples}
The $\Span$ construction was analyzed further in \cite{Str00}.
In particular, for any small category $\ca C$ in place of $\G$,
there is a $2$-adjunction
\[ \xymatrix{{\CAT} \ar@/_{1pc}/[rr]_-{\Span_{\ca C}}
\save \POS?="bot" \restore
&& {[\op {\ca C},\CAT]} \ar@/_{1pc}/[ll]_-{\EL}
\save \POS?="top" \restore
\POS "bot"; "top" **@{}; ?*{\perp}} \]
where $\Span_{\ca C}(\ca E)(C) = [\op {({\ca C}/C)},\ca E]$,
and $\EL(X)$ is the following category:
\begin{itemize}
\item  objects are pairs $(C,x)$ where $C \in \ca C$
and $x \in X(C)$.
\item  morphisms $(C,x)\rightarrow(D,y)$
are pairs $(f,\alpha)$ where $f:D{\rightarrow}C$ in $\ca C$,
and $\alpha:X(f)(x){\rightarrow}y$ in $X(D)$.
\item  compositions and identities are inherited in the obvious way
from $\ca C$ and the categories $X(C)$.
\end{itemize}
When $X$ is discrete, that is, as a functor factors through $\SET$,
then $\EL(X)=\op {\el(X)}$, the dual of the usual category of elements of $X$.
If moreover, $X$ is small, that is, factors through $\Set$,
and $\ca E = \Set$, then we have
\begin{eqnarray*}
[\op {\ca C},\CAT](X,\Span_{\ca C}(\Set))
& \iso & \CAT(\EL(X),\Set) \\
& = & \CAT(\op {\el(X)},\Set) \\
& \catequiv & {\PSh {\ca C}}/X
\end{eqnarray*}
this last step being a well known equivalence of categories, pseudo natural in $X$.
Let $(M,\eta,\mu)$ be a cartesian monad on $\PSh {\ca C}$,
and recall the category $\Coll M$ from \cite{Kel92} and \cite{WebGen},
which is the full subcategory
of $[\PSh {\ca C},\PSh {\ca C}]/M$ consisting of the cartesian
natural transformations.
Recall also that $\Coll M$ has a strict monoidal structure:
\begin{itemize}
\item  $\eta:1{\rightarrow}M$ is the unit.
\item  $\phi\tensor\psi$ is the composite
\[ \xymatrix{ST \ar[r]^-{\phi\psi} & MM \ar[r]^{\mu} & M}, \]
\end{itemize}
and that evaluation at $1$ provides an equivalence
of categories ${\PSh {\ca C}}/M(1) \catequiv \Coll M$.
That is, given a cartesian monad $(M,\eta,\mu)$
on $\PSh {\ca C}$, we have
\[ [\op {\ca C},\CAT](M(1),\Span_{\ca C}(\Set))
\catequiv \Coll M \]
We now present the higher operads of \cite{Bat98}.
\begin{examples}\label{ex:hops}
\begin{enumerate}
\item  A $\ca T$-operad $(p,\iota,\sigma)$ in $A$
amounts to a higher operad in an augmented monoidal
globular category $A$ in the sense of \cite{Bat98},
subject to one caveat.
That is, the above definition
is in fact more general than that presented in \cite{Bat98}.
The difference is that in \cite{Bat98},
further hypotheses on $A$ are required, namely,
that $A$ has globular coproducts
which are compatible with
the monoidal pseudo $\ca T$-algebra structure of $A$
(see \cite{Bat98} for further elaboration).
In \cite{AnOp}, these hypotheses are seen as another instance
of a general notion of distributive monoidal pseudo algebra.
These further hypotheses induce a monoidal structure
on the category
$\GLOBCAT(\ca T(1),A)$ (\cite{Bat98} Theorem 6.1)
and operads were defined by Batanin to be monoids
in this monoidal category.
It can be verified directly
that the category of monoids in $\GLOBCAT(\ca T(1),A)$
is isomorphic to $\Op T A$.
\label{hops-gen}
\item  For the case $A = \Span(\Set)$ of (\ref{hops-gen}),
we shall continue to regard $\ca T$ as a $2$-monad
on $\GLOBCAT$, and $\ca T_0$ as a monad on $\Glob$.
Now $\ca T(1)=\ca T_0(1)$, and the equivalence
\[ [\op {\G},\CAT](\ca T_0(1),\Span(\Set))
\catequiv \Coll {\ca T_0} \]
is in fact a monoidal equivalence.
Thus, a $\ca T$-operad in $\Span(\Set)$
amounts to a cartesian monad morphism $\phi_0:R_0{\rightarrow}{\ca T_0}$,
and algebras for this operad amount to algebras for the monad $R_0$.
We shall call such an operad $\phi_0$ a \emph{basic} higher operad.
There is a basic higher operad whose algebras are weak $\omega$-categories.
\label{basic}
\item  By (\ref{ex:mgc})(\ref{symmon})
one can consider $\ca T_{(n)}$-operads within symmetric monoidal categories,
where $n{\geq}3$.
Such examples are important for the applications of higher operads
to the study of loop spaces, see \cite{Bat02} and \cite{Bat03}.
\label{sym}
\item  Let $\phi_0:R_0{\rightarrow}{\ca T_0}$ be a basic higher operad.
Applying the $2$-functor $\Cat$ (take category objects),
and shifting up to the next set-theoretic universe,
we have cartesian $2$-monad morphism $\phi:R{\rightarrow}{\ca T}$.
The induced forgetful $2$-functors
${\Alg {\ca T}}{\rightarrow}{\Alg R}$,
${\PsAlg {\ca T}}{\rightarrow}{\PsAlg R}$, and
${\NormPsAlg {\ca T}}{\rightarrow}{\NormPsAlg R}$
preserve products, and so in particular $\phi$ induces a forgetful $2$-functor
\[ {\PsMon {\NormPsAlg {\ca T}}}{\rightarrow}{\PsMon {\NormPsAlg R}} \]
ensuring a ready supply of examples of monoidal pseudo $R$-algebras.
So it is potentially interesting to
consider $R$-operads. For example, one could consider the case
where $R_0$ is the weak $\omega$-category monad,
for a weakened version of higher operad.
\label{spec}
\end{enumerate}
\end{examples}
%

\section{Symmetric variants of higher operads}\label{higher-sym-ops}

In this section, the symmetric analogues of Batanin's higher operads
are described.
In order to do so, new examples of $2$-monads
on $\GLOBCAT$ are described, which
blend together the $2$-monad $\ca T$,
with an appropriate $2$-monad $\ca C$ on $\CAT$.
For instance taking $\ca C$
to be $\ca B$, the braided monoidal category $2$-monad,
the blend alluded to here mixes the combinatorics
of trees and pasting diagrams encapsulated by $\TREE$,
with that of braids, and the operad notion
corresponding to this new monad is a braided analogue
of higher operad.
This construction hinges on two things:
\begin{enumerate}
\item  The underlying $2$-functor of $C$
preserves pullbacks. This is easily observed
directly for the examples of interest: $\ca M$, $\ca B$ and $\ca S$.
\item  An alternative description of $\wCat$,
and more generally $\Alg R_0$
for a basic higher operad $\phi_0:R_0{\rightarrow}\ca T_0$
(as in (\ref{ex:hops})(\ref{basic})),
as models for a finite connected limit sketch.
\end{enumerate}
This last point is not particularly surprising,
at least for strict $\omega$-categories.
Already from \cite{Str86},
one knows that $\wCat$ is the category of $\Set$-valued
models for some sketch.
However, from the work of Clemens Berger \cite{Ber00},
as we shall now explain,
one can obtain this sketch directly from the monad $\ca T_0$.
Moreover, this procedure generalises to any basic higher operad.

Following \cite{Ber00} we regard $\Theta_0$ as a
Grothendieck site by taking
covering families to be jointly epimorphic families
of morphisms.
Denote by $\Shv {\Theta_0}$ the category of sheaves
on the site $\Theta_0$.
Let $\phi_0:R_0{\rightarrow}{\ca T_0}$ be a basic higher operad,
and denote by $\Theta_R$ the full subcategory of $\Alg {R_0}$
whose objects are the globular cardinals,
and write $i_R:{\Theta_R}{\hookrightarrow}{\Alg {R_0}}$
for the inclusion.
Via the left adjoint $\Glob{\rightarrow}{\Alg {R_0}}$
to the forgetful functor,
one can identify $\Theta_0$ as a subcategory of $\Theta_R$.
Since $R_0$ is a finitary monad on $\Glob$,
$\Alg R_0$ is locally finitely presentable and so cocomplete.
Thus one obtains a ``hom-tensor'' adjunction
\[ \xymatrix{{[\op \Theta_R,\Set]} \ar@/^1pc/[rr]^-{\ca L_R}
\save \POS?="top" \restore
&& {\Alg {R_0}} \ar@/^1pc/[ll]^-{\ca N_R}
\save \POS?="bot" \restore
\POS "top";"bot" **@{}; ?*{\perp}} \]
where $\ca L_R$ is the left kan extension of $i_R$ along
the yoneda embedding,
and $\ca N_R(X)(T) = {\Alg {R_0}}(i_R(T),X)$.
\begin{definition}\cite{Ber00}
A $\Theta_R$-model is a presheaf $F \in [\op \Theta_R,\Set]$
whose restriction to $\Theta_0$ is a sheaf.
Denote by $\Mod {\Theta_R}$ the full subcategory
of $[\op \Theta_R,\Set]$ consisting of the $\Theta_R$-models.
\end{definition}
\begin{theorem}\cite{Ber00}\label{thm:tmod}
For any basic higher operad
$\phi_0:R_0{\rightarrow}{\ca T_0}$:
\begin{enumerate}
\item  $\ca N_R$ is fully faithful.
\item  The adjunction $\ca L_R \ladj \ca N_R$
restricts to an equivalence
$\Mod {\Theta_R} \catequiv \Alg {R_0}$.
\label{mod-alg-equiv}
\end{enumerate}
\end{theorem}
Note that in general, the fully faithfulness of $\ca N_R$
is equivalent to the density of $i_R$.
\begin{examples}
\begin{enumerate}
\item  For the basic higher operad $\eta:1{\rightarrow}{\ca T}$,
(\ref{thm:tmod})(\ref{mod-alg-equiv})
gives an equivalence $\Glob \catequiv \Shv {\Theta_0}$.
This equivalence can also be seen as a basic consequence
of the Giraud theorem from topos theory (see \cite{MM} pg 589),
since $\Theta_0$, which contains the representables, generates $\Glob$.
\item  $\Theta_{\ca T}$ is denoted as $\Theta$ in \cite{Ber00},
where it was shown to coincide with Joyal's category $\Theta$
from \cite{Joy97}.
\end{enumerate}
\end{examples}
For our purposes a mild variation on this characterisation
of the categories $\Alg {R_0}$ is necessary, namely, as the models
of a finite connected limit sketch.
First we recall the definition of limit sketch
and models thereof.
\begin{definition}
A limit sketch is a $4$-tuple $\ca D = (D,I,F,c)$
where $D$ is a category, $I$ is a set,
$F$ is an $I$-indexed set of functors
$F_i:J_i{\rightarrow}D$,
and $c$ is an $I$ indexed set of cones
$c_i:\Delta(x_i){\implies}F_i$
(where $\Delta(x_i)$ denotes the functor constant at $x_i$).
We call the set $F$ the \emph{diagrams},
and the set $c$ the \emph{distinguished cones}
for the sketch $\ca D$.
Let $\ca E$ be a category with limits of functors out of $J_i$.
The category $\ArbMod {\ca D} {\ca E}$,
of $\ca E$-valued models of $\ca D$,
is the full subcategory of $[D,\ca E]$
consisting of the functors $D{\rightarrow}\ca E$
which take the cones $c_i$ to limiting cones.
Denote by $\Mod {\ca D}$ the category $\ArbMod {\ca D} \Set$.
\end{definition}
\begin{examples}\label{sk-ex}
\begin{enumerate}
\item  It is well known that limit sketches subsume Grothendieck topologies,
for, let $D$ be a category, and $\ca J$ a Grothendieck topology on $D$.
Note that for each sieve $\alpha{\hookrightarrow}D(-,x) \in {\ca J}$,
one gets a diagram $\el(\alpha){\rightarrow}D$ as the discrete fibration
corresponding to $\alpha$,
and a cocone $c$ for this diagram with components $c_{(y,f)}=\alpha_y(f)$.
In this way one gets a distinguished \emph{cone} in $\op D$ for each
sieve in $\ca J$,
and so a limit sketch whose underlying category is $\op D$.
By definition, $\Set$-valued models for this sketch are sheaves
for the Grothendieck topology $\ca J$.
\label{gtop->sk}
\item  Let $\phi_0:R_0{\rightarrow}{\ca T_0}$ be a basic higher operad.
By the above example,
one has a limit sketch whose underlying category is $\op \Theta_0$
from the Grothendieck topology described above (covering maps are
jointly epimorphic families).
By composing with the inclusion ${\Theta_0}{\rightarrow}{\Theta_R}$,
one has a limit sketch whose underlying category is $\op \Theta_R$,
and by definition, $\Mod {\Theta_R}$ is the category of $\Set$-valued
models for this sketch.
This sketch does not typically arise from a Grothendieck topology.\footnotemark{
\footnotetext{For example, $\Alg {\ca T_0} = \wCat$ is not a Grothendieck topos,
in fact, one can show that it is not even a regular category.}}
We shall abuse notation and refer to this sketch as $\Theta_R$,
even though the underlying category of this sketch is $\op \Theta_R$.
\label{th->sk}
\end{enumerate}
\end{examples}
\begin{definition}
A limit sketch $\ca D = (D,I,F,c)$, is a connected limit
sketch when the categories $J_i$ (that is, the domains
of the $F_i$) are connected.
$\ca D$ is a finite limit sketch when the $J_i$ have a finite initial
subcategory.
\end{definition}
For a finite connected limit sketch,
the distinguished limiting cones may be regarded as iterated pullbacks.
More precisely, for such a sketch $\ca D$,
one can define $\ArbMod {\ca D} {\ca E}$
as long as $\ca E$ has pullbacks,
and composition with a pullback preserving
functor ${\ca E}{\rightarrow}{\ca E'}$
induces ${\ArbMod {\ca D} {\ca E}}{\rightarrow}{\ArbMod {\ca D} {\ca E'}}$.
However, in the general context of (\ref{sk-ex})(\ref{gtop->sk}),
there is nothing forcing the diagram
corresponding to
an arbitrary sieve $\alpha{\hookrightarrow}D(-,x) \in {\ca J}$,
to be finite or connected. We shall now show
that for $\Theta_0$ with the given Grothendieck topology,
that this is indeed the case, and so, (\ref{sk-ex})(\ref{th->sk})
is a finite connected limit sketch.

For a linearly ordered set $X$, and $x \in X$,
we shall write $x^+$ for the successor of $x$,
which exists as long as $x$ is not the maximum element of $X$.
\begin{lemma}\label{char-succ}
Let $X$ be a globular cardinal.
Regarding $x \in X_n$ as an element of $\Sol(X)$,
and assuming $x$ is not the maximum element,
we have
\[ x^+ = \left\{
\begin{array}{rcl}
y && \textnormal{if $s(y)=x$} \\
t(x) && \textnormal{otherwise}
\end{array}\right. \]
\end{lemma}
\begin{proof}
Suppose that $x=s(y)$ and $x{\solleq}z{\solleq}y$.
If $z{\neq}x$ then we have $x{\prec}a{\solleq}z$,
and so $a$ must be either $y$ or $t(x)$.
In the first case, $a=y$, we have $y{\solleq}z{\solleq}y$
and so $y=z$. On the other hand, if $a=t(x)$, note that
$t(x)=ts(y)=tt(y)$, and
$t(x){\solleq}z{\solleq}y{\solleq}t(y){\solleq}tt(y)$,
so that $y=z$ also.
Thus if $x=s(y)$, we have $x^+=y$.
On the other hand suppose that there is no $y$
such that $x=s(y)$.
If there were no $t(x)$, then $x{\in}X_0$, and
$X$ would be the globular set with one $0$-cell
and no other cells, in which case $x$ is the maximum.
Now, suppose $x{\solleq}z{\solleq}t(x)$
and $x{\neq}z$.
Then we have $x{\prec}a{\solleq}z$ and $a$ is forced to be $t(x)$.
Thus, $t(x){\solleq}z{\solleq}t(x)$ and so $z=t(x)$.
Thus if there is no $y$ such that $x=s(y)$,
then $x^+=t(x)$.
\end{proof}
\begin{corollary}\label{succ-pres}
Let $f:X{\rightarrow}Y$ in $\Theta_0$.
Then $\Sol(f)(x^+)=\Sol(f)(x)^{+}$
for all non-maximal elements $x$.
\end{corollary}
\begin{proof}
By (\ref{char-succ}) the successor operation
for $\Sol(X)$ is expressed in terms of the sources and targets for $X$,
which are preserved by $f$ since it is a morphism
of globular sets.
\end{proof}
Given non-empty finite linear orders $X$, $Y$ and $Z$,
and successor-preserving maps
\[ \xymatrix{X \ar[r]^-{f} & Y & Z \ar[l]_-{g}}, \]
the pullback in $\PreOrd$ of these maps
is a finite linear order. It will simply be formed
as the intersection of the images of $f$ and $g$.
Recall that in any category $\ca E$ with pullbacks and an initial object,
arrows $f$ and $g$ as above are said to be \emph{disjoint}
when their pullback is the initial object of $\ca E$.
\begin{proposition}\label{pb-glob-card}
$\Theta_0$ has pullbacks of pairs of maps
which are non-disjoint in $\Glob$.
\end{proposition}
\begin{proof}
Let
\[ \xymatrix{P \ar[r] \ar[d] & X \ar[d]^{f} \\ Y \ar[r]_-{g} & Z} \]
be a pullback square in $\Glob$, $f$ and $g$ non-disjoint,
and X, Y and Z globular cardinals.
Applying $\Sol$, which preserves pullbacks and initial objects,
to this pullback square, exhibits $\Sol(P)$ as a pullback
of non-disjoint successor-preserving maps between finite linear orders.
Thus, $\Sol(P)$ is a non-empty finite linear order, and so $P$ is a globular
cardinal.
\end{proof}
\begin{proposition}\label{thsk-fin-conn}
The limit sketch of (\ref{sk-ex})(\ref{th->sk})
is a finite connected limit sketch.
\end{proposition}
\begin{proof}
Let $F$ be a jointly epimorphic family of maps in $\Theta_0$ with codomain $X$,
let ${\alpha}{\hookrightarrow}\Theta_0(-,X)$ be the sieve
generated by $F$.
We must show that $\el(\alpha)$
is connected and has a finite final subcategory.
First note that $X$ is non-empty since it is a globular cardinal, and so
$F$ and $\alpha$ are non-empty also.
Let $f$ and $f'$ be a pair of maps in $F$.
Then there will be a finite sequence $(f_0, ..., f_n)$
of maps from $F$, such that all consecutive pairs of maps
in the sequence $(f,f_0, ..., f_n, f')$ are non-disjoint,
since $X$ is finite and $F$ is a jointly epimorphic family.
By (\ref{pb-glob-card}), one can take the joint pullback
of the maps $(f,f_0, ..., f_n, f')$ in $\Theta_0$
to exhibit $\el(\alpha)$ as connected.
Since all maps in $F$ are monic by (\ref{theta0-old}),
and $X$ has only finitely many subobjects,
$\el(\alpha)$ contains only finitely many maps
up to isomorphism in $\Glob/X$.
That is, $\el(\alpha)$ is actually equivalent to a finite category.
\end{proof}
\begin{corollary}\label{cor:cat-alg-basic}
Let $\phi_0:R_0{\rightarrow}{\ca T_0}$ be a basic higher operad.
Then
\[ \Alg R \catequiv \ArbMod {\Theta_R} \Cat \]
\end{corollary}
\begin{proof}
By (\ref{thsk-fin-conn}), we can apply the $2$-functor
${\Cat}:{\CART}{\rightarrow}{\TWOCAT}$,
which takes category objects (see (\ref{Monad-ex})(\ref{Cat-of-monad})),
to the equivalence
$\Mod {\Theta_R} \catequiv \Alg {R_0}$
of (\ref{thm:tmod})(\ref{mod-alg-equiv}).
Clearly, $\ArbMod {\Theta_R} \Cat \iso \Cat(\Mod {\Theta_R})$.
\end{proof}
\begin{remark}
By (\ref{thsk-fin-conn})
it makes sense to take models of any basic operad
$\phi:R_0{\rightarrow}{\ca T_0}$
in any category $\ca E$ with pullbacks.
For example,
in this way one can speak about weak $\omega$-categories
internal to $\ca E$.
\end{remark}
Any $2$-monad $(C,\eta,\mu)$ on $\Cat$
may be regarded as a $2$-monad $(C_{\G},\eta_{\G},\mu_{\G})$
on $[{\op \G},\Cat]$
by composition, that is, the components of $\eta_{\G}$ and $\mu_{\G}$
for $X \in [{\op \G},\Cat]$ are:
\[ \xymatrix{{\op \G} \ar[rr]^-{X} &&
{\Cat} \ar@/^{2pc}/[rr]^-{1} \save \POS?="domunit" \restore
\ar[rr]|-{C}
\save \POS?="codunit" \restore \save \POS?="codmult" \restore
\ar@/_{2pc}/[rr]_-{CC} \save \POS?="dommult" \restore && {\Cat}
\POS "domunit";"codunit" **@{} ?(.3) \ar@{=>}^{\eta} ?(.7)
\POS "dommult";"codmult" **@{} ?(.3) \ar@{=>}^{\mu} ?(.7)}, \]
and we shall see that this $2$-monad distributes with $\ca T$
whenever $C$ preserves pullbacks.
First, we shall clarify what we mean by a distributive law between $2$-monads,
since there are various notions that one could use.

Recall that a distributive law between monads $S$ and $T$,
is a natural transformation
${\lambda}:TS{\rightarrow}ST$, satisfying some axioms,
which enable one to define a monad structure on the composite $ST$.
This is done is such a way that algebra structures for $ST$,
amount to compatible $S$ algebra and $T$ algebra structures.
In \cite{Str72} it was shown that
the theory of monad distributive laws
can be developed internal to a $2$-category.
In particular, when the $2$-category $\ca K$ in question has
certain weighted limits called Eilenberg-Moore objects,
distributive laws ${\lambda}:TS{\rightarrow}ST$
between monads $S$ and $T$ on $A \in \ca K$,
correspond to liftings of the monad $S$
to the Eilenberg-Moore object $A^T$ (object
of $T$-algebras).
In more detail, recall that the Eilenberg Moore object
includes a ``forgetful one-cell''
$u:A^T{\rightarrow}A$. A lifting of $f:A{\rightarrow}A$
is an $\overline{f}$ making
\[ \xymatrix{{A^T} \ar[r]^-{\overline{f}} \ar[d]_{u} & {A^T} \ar[d]^{u} \\
A \ar[r]_-{f} & A} \]
commute.
Given liftings $\overline{f_1}$ and $\overline{f_2}$
of $f_1$ and $f_2$, and a $2$-cell $\phi:f_1{\implies}f_2$,
a lifting of $\phi$ from $\overline{f_1}$ to $\overline{f_2}$
is a $2$-cell $\overline{\phi}$ making
\[ \xymatrix{{A^T} \ar@/^1pc/[r]^-{\overline{f_1}} \save \POS?="dompb" \restore
\ar@/_1pc/[r]_-{\overline{f_2}} \save \POS?="codpb" \restore
\ar[dd]_{u} & {A^T} \ar[dd]^{u} \\ \\
A \ar@/^1pc/[r]^-{f_1} \save \POS?="domp" \restore
\ar@/_1pc/[r]_-{f_2} \save \POS?="codp" \restore & A
\POS "dompb"; "codpb" **@{} ?(.3) \ar@{=>}^{\overline{\phi}} ?(.7)
\POS "domp"; "codp" **@{} ?(.3) \ar@{=>}^{\phi} ?(.7)} \]
commute. So to give a distributive law $TS{\rightarrow}ST$,
is to give a lifting in this sense, of
all the data of the monad $(S,\eta,\mu)$ on $A$,
to a monad $(\overline{S},\overline{\eta},\overline{\mu})$
on $A^T$.
See \cite{Str72} and \cite{LS00} for further elaboration.
For us, the $2$-category $\ca K$ is that of $\CAT$-enriched categories:
monads in $\ca K$ are $2$-monads, and the Eilenberg-Moore
object of $T$ is the $2$-category $\sAlg T$ of strict $T$-algebras,
strict algebra morphisms, and algebra $2$-cells.
When there is a distributive law $TS{\rightarrow}ST$,
in this sense between $2$-monads $S$ and $T$,
we shall say that $S$ \emph{distributes} with $T$.
\begin{theorem}\label{thm:globdist}
Let $(C,\eta,\mu)$ be a $2$-monad on $\Cat$
such that $C$ preserves pullbacks,
and $\phi_0:R_0{\rightarrow}{\ca T}_0$
be a basic higher operad.
Then the $2$-monad $C_{\G}$ distributes with $R$.
\end{theorem}
\begin{proof}
By the theory of distributive laws
it suffices to exhibit a lifting of $(C_{\G},\eta_{\G},\mu_{\G})$
to $\sAlg R$, and by (\ref{cor:cat-alg-basic})
we have the equivalence $\sAlg R \catequiv \ArbMod {\Theta_R} \CAT$.
Since $C$ preserves pullbacks, the $2$-monad on $[\op \Theta_R,\CAT]$
with components
\[ \xymatrix{{\op \Theta_R} \ar[rr]^-{X} &&
{\CAT} \ar@/^{2pc}/[rr]^-{1} \save \POS?="domunit" \restore
\ar[rr]|-{C}
\save \POS?="codunit" \restore \save \POS?="codmult" \restore
\ar@/_{2pc}/[rr]_-{CC} \save \POS?="dommult" \restore && {\CAT}
\POS "domunit";"codunit" **@{} ?(.3) \ar@{=>}^{\eta} ?(.7)
\POS "dommult";"codmult" **@{} ?(.3) \ar@{=>}^{\mu} ?(.7)}, \]
restricts to $\ArbMod {\Theta_R} \CAT$,
and is by definition a lifting of $(C_{\G},\eta_{\G},\mu_{\G})$.
\end{proof}
\begin{example}
Let $\ca E$ be a category with finite limits,
and $\phi_0:R_0{\rightarrow}{\ca T_0}$ be a basic higher operad.
Then by (\ref{ex:hops})(\ref{spec}) and
(\ref{ex:mgc})(\ref{span}),
$\Span(\ca E)$ has a monoidal pseudo $R$-algebra structure,
with the additional pseudo monoid structure
being given dimensionwise by pointwise cartesian product.
Thus, this additional pseudo monoid structure
may be regarded as a compatible pseudo $\ca S_{\G}$ algebra structure.
In this way, $\Span(\ca E)$ is canonically a
pseudo $\ca S_{\G}R$-algebra.
As in (\ref{ex:mon-psalg})(\ref{symmon}),
the pseudo monoid part of a monoidal pseudo $\ca S_{\G}R$-algebra
encodes no new information.
Thus, for any basic higher operad $\phi_0:R_0{\rightarrow}{\ca T_0}$,
and category $\ca E$ with finite limits,
$\Span(\ca E)$ is canonically a
monoidal pseudo $\ca S_{\G}R$-algebra.
Via the forgetful functors induced by
the obvious monad morphisms $\ca M{\rightarrow}{\ca B}{\rightarrow}{\ca S}$,
$\Span(\ca E)$ may be regarded as a monoidal pseudo algebra
also for the monads $\ca M_{\G}R$ and $\ca B_{\G}R$.
\end{example}
By (\ref{thm:globdist}) we can consider
$\ca M_{\G}\ca T$-operads, $\ca B_{\G}\ca T$-operads,
and $\ca S_{\G}{\ca T}$-operads
-- natural higher globular analogues
of non-symmetric operads, braided operads and symmetric operads respectively.
For that matter, one may replace $\ca T$ by $R$, for an arbitrary
basic higher operad $\phi_0:R_0{\rightarrow}{\ca T_0}$.
Thus, any higher dimensional categorical structure, which is describable
by a basic higher operad,
automatically comes equipped with its own analogous notions
of non-symmetric, braided, and symmetric operad.

\section{ACKNOWLEDGMENTS}

This work was done while the author was employed as a
post doctoral research fellow
at the Department of Mathematics and Statistics
in the University of Ottawa.
I am thankful to the members of this department
for providing a friendly and stimulating work environment.

I am grateful to Michael Batanin, Ross Street,
Steve Lack and Claudio Hermida for all
that they have taught me about higher dimensional category theory
and related subjects.
This work would not have evolved at all were it not for them.
Here in Ottawa, I would like to thank Rick Blute for organising
the informal homotopy-logic seminars.
It was in trying to describe higher operads in these meetings
that I happened upon the simpler and more general formalism presented
in this paper.

\bibliography{highercats}
\end{document}